\documentclass[11pt,reqno]{amsart}

\theoremstyle{plain}
\newtheorem{theorem} {Theorem}[section]
\newtheorem{lemma}[theorem] {Lemma}
\newtheorem{proposition}[theorem] {Proposition}
\newtheorem{corollary}[theorem] {Corollary}

\theoremstyle{definition}

\newtheorem{example} [theorem]{Example}

\theoremstyle{remark}
\newtheorem{remark}[theorem] {Remark}

\numberwithin{equation}{section}

\usepackage{amsfonts}
\usepackage{amssymb}
\usepackage{amscd}

\newcommand{\R}{{\mathbb R}}
\newcommand{\Z}{{\mathbb Z}}
\newcommand{\N}{{\mathbb N}}

\newcommand{\PP}{{\mathcal P}}

\newcommand{\C}{{\mathcal C}}
\newcommand{\CC}{{\mathbb C}}

\newcommand{\E}{{\mathcal E}}
\newcommand{\B}{{\mathcal B}}
\newcommand{\BB}{{\mathfrak B}}
\newcommand{\al}{{\alpha}}
\newcommand{\la}{{\lambda}}
\newcommand{\sa}{{\sigma}}

\newcommand{\iy}{{\infty}}
\newcommand{\vphi}{{\varphi}}
\newcommand{\vep}{{\varepsilon}}
\newcommand{\g}{{\gamma}}
\newcommand{\de}{{\delta}}
\newcommand{\Om}{{\Omega_m}}

\newcommand{\be}{{\beta}}

\newcommand{\bna}{\begin{eqnarray}}
\newcommand{\ena}{\end{eqnarray}}
\newcommand{\ba}{\begin{eqnarray*}}
\newcommand{\ea}{\end{eqnarray*}}
\newcommand{\beq}{\begin{equation}}
\newcommand{\eeq}{\end{equation}}

\begin{document}

\title[Constants in Multivariate Inequalities]
{Sharp Constants of Approximation Theory. II. Invariance Theorems and Certain Multivariate Inequalities
of Different Metrics}
\author{Michael I. Ganzburg}
 \address{Department of Mathematics\\ Hampton University\\ Hampton,
 VA 23668\\USA}
 \email{michael.ganzburg@hamptonu.edu}
 \keywords{Sharp constants, multivariate Markov-Bernstein-Nikolskii
  type inequality, algebraic polynomials,
 entire functions of exponential type, weighted spaces.}
 \subjclass[2010]{Primary 41A17, 41A63, Secondary 26D10}

 \begin{abstract}
We prove invariance theorems for general inequalities
 of different metrics and apply them to
limit relations between the sharp constants in the
multivariate
Markov-Bernstein-Nikolskii
type inequalities
 with the polyharmonic operator for algebraic polynomials on the unit sphere
 and the unit ball in $\R^m$ and the corresponding constants for
 entire functions of spherical type on $\R^m$.
Certain relations in the univariate weighted spaces are discussed as well.
 \end{abstract}
 \maketitle

 \section{Introduction}\label{S1}
\setcounter{equation}{0}
\noindent
We continue the study of the sharp constants in multivariate inequalities
of approximation theory
that began in \cite{G2018}. In this paper we prove invariance theorems for multivariate inequalities
 of different metrics and apply them to
limit relations between the sharp constants in the
multivariate
 Markov-Bernstein-Nikolskii
type inequalities for algebraic polynomials and entire functions of exponential type.
In addition, we discuss the asymptotic behavior of certain sharp constants
 in univariate weighted spaces.\vspace{.12in}\\
\textbf{Notation.}
Let $\R^m,\,m\ge 1,$ be the Euclidean $m$-dimensional space with elements
$x=(x_1,\ldots,x_m),\linebreak y=(y_1,\ldots,y_m)$,
 the inner product $(x,y):=\sum_{j=1}^mx_jy_j$,
and the norm $\vert x\vert:=\sqrt{(x, x)}$.
Next, $\CC^m:=\R^m+i\R^m$ is the $m$-dimensional complex space with elements
$z=x+iy=(z_1,\ldots, z_m)$
and the norm $\vert z\vert:=
\sqrt{\sum_{j=1}^m\vert z_j\vert^2}=
\sqrt{\vert x\vert^2+\vert y\vert^2}$;
$\Z^m_+$ denotes the set of all integral lattice points
in $\R^m$ with nonnegative coordinates.
In addition, we use a multi-index $\al=(\al_1,\ldots,\al_m)\in \Z^m_+$  with
 $\vert\al\vert:=\sum_{j=1}^m\al_j$
  and $x^\al:=x_1^{\al_1}\cdot\cdot\cdot x_m^{\al_m}$.
We also use the Frobenius norm
$\|A\|_F:=\left(\sum_{k,j=1}^m a_{k,j}^2\right)^{1/2}$
of an $m\times m$ matrix $A=\left[a_{k,j}\right]_{k,j=1}^m$
with real elements.
Given $M>0$, let
$\BB^m(M):=\{x\in\R^m: \vert x\vert\le M\},\,\BB^m:=\BB^m(1)$, and
$S^{m-1}:=\{x\in\R^m: \vert x\vert=1\}$
be the $m$-dimensional ball of radius $M$,
the unit  $m$-dimensional ball,
 and the unit $(m-1)$-dimensional sphere in $\R^m$, respectively.
Next, let $\vert E\vert_k$ denote the $k$-dimensional Lebesgue measure of a
$k$-dimensional measurable set $E\subseteq\R^m$. In particular,
$\vert S^{m-1}\vert_{m-1}=2\pi^{m/2}/\Gamma(m/2),\,m\ge 2$, and
$\vert S^0\vert_0:=2$.
In addition, we use generic notation $\lfloor a \rfloor,\,\Gamma(z)$
 and $\mbox{B}(\al,\be)$ for the floor function, the gamma function
 and the beta function, respectively.

Throughout the paper $C,\,C_0,\,C_1,\ldots$ denote positive constants independent
of essential parameters.
 Occasionally we indicate dependence on certain parameters. The same symbol $C$ does not
 necessarily denote the same constant in different occurrences.
  \vspace{.12in}\\
\textbf{ Markov-Bernstein-Nikolskii Type Inequalities.}
Limit relations between sharp constants
 in the univariate  Markov-Bernstein-Nikolskii type inequalities for trigonometric
 and algebraic polynomials and entire functions of exponential type
 were studied by Taikov \cite{T1965, T1993}, Gorbachev \cite{G2005},
 Levin and Lubinsky \cite{LL2015a, LL2015b}, the author and Tikhonov \cite{GT2017},
 and the author \cite{G2017}.
Detailed surveys of the univariate  Markov-Bernstein-Nikolskii type inequalities
for trigonometric and algebraic polynomials and entire functions of exponential type
were presented in \cite{GT2017, G2017}. The corresponding multivariate problems
 were recently studied by  Dai, Gorbachev, and Tikhonov
 \cite{DGT2018} and the author \cite{G2018}.
 The following publications \cite{K1974, ABD2018, Gor2018} are also closely
  related to the subject of the paper.

 The purpose of this paper is threefold.
 First, we extend invariance theorems of approximation theory,
 proved by the author and Pichugov \cite{GP1981} and by the author
 \cite{G2008}, to the generalized
  Markov-Bernstein-Nikolskii type inequalities.
 These results are presented and proved in Section
 \ref{S2} (see Theorems \ref{T2.1}
 and \ref{T2.2}).
 In particular, invariance theorems reduce certain
 multivariate inequalities to univariate ones in weighted metrics.
 Certain special cases are discussed in Section
 \ref{S3} (see Examples \ref{Ex2.10}, \ref{Ex2.12},
  \ref{Ex2.14} and Corollaries \ref{C2.11}, \ref{C2.13}, \ref{C2.15}).

 Second, in Section \ref{S4} we obtain limit relations between the sharp constants
 in the multivariate  Markov-Bernstein-Nikolskii type inequalities
 with the polyharmonic operator for polynomials on the unit sphere
 and the unit ball in $\R^m$ and the corresponding constants for
 entire functions of spherical type on $\R^m$ (see
 Corollaries \ref{C4.4} and \ref{C4.5}).
 In particular, the major part of Corollary \ref{C4.4}
 is the following statement:
 if $\Delta$ is the Laplace operator, then the minimal (sharp) constant
 $C_n^*$ in the inequality
 \beq\label{E1.1}
 \vert \Delta^NP(0)\vert\le C_n^*
 \left(\int_{\vert x\vert\le 1}\vert P(x)\vert^pdx\right)^{1/p}
 \eeq
 over all polynomials $P$ of degree at most $n$ in $m$ variables
 and the minimal constant
 $C^*$ in the inequality
 \ba
 \| \Delta^Nf\|_{L_\iy(\R^m)}\le C^*
 \left(\int_{\R^m}\vert f(x)\vert^pdx\right)^{1/p}
 \ea
 ($p\in[1,\iy],\,N=0,\,1\,\ldots$) over all entire functions $f$
 of spherical type $1$ are connected by the limit relation
 $\lim_{n\to\iy}n^{-2N-m/p}C_n^*=C^*.$

 A similar result on the unit sphere is proved in Corollary \ref{C4.5}.
 A special case of Corollary \ref{C4.5} for $N=0$ was established in
 \cite{DGT2018}.

 The proofs of multivariate results
 are based on the invariance theorems and certain univariate results.
By using the invariance theorems, the limit multivariate relations
  can be reduced to the relations between sharp constants in the
  univariate  Markov-Bernstein-Nikolskii type inequalities for
  algebraic polynomials
  with the Bessel and Gegenbauer differential operators in weighted
  $L_p$-spaces on $[-1,1],\,1\le p\le\iy$, and the corresponding constants
  for univariate entire functions of exponential type.

  Third, limit  relations between the univariate sharp constants
  in more general weighted spaces
 are presented in Theorems \ref{T4.1} and \ref{T4.3}.
 Their proofs are given in Section \ref{S6}.
    Note that the proofs are based on an approach
    to limit relations between sharp constants
    developed in \cite{GT2017, G2017, G2018}.

 Surprisingly, a special case of limit relations between
 the univariate sharp constants in weighted spaces is the asymptotic
 relation for the sharp constant in the classical inequality for  univariate
 polynomials of different metrics (see Corollary \ref{C4.6}; cf.
 \cite[Theorem 1.4 and p. 94]{G2017}).

Special cases of the results in Section \ref{S4}
 are obtained earlier in \cite{LL2015a, G2017, DGT2018}.
Section \ref{S5} contains certain properties of
 entire functions of exponential type and polynomials
 that are needed for the proofs.

\section{General Invariance Theorems}\label{S2}
 \noindent
\setcounter{equation}{0}
The averaging (or symmetrization) principle is well known
in approximation theory (see, e.g., \cite{L1966, Ch1969, A1996,
GP1981, G2008} and references therein).
In particular, general and special invariance theorems for the error
of best approximation were proved in \cite{GP1981, G2008}.
In this section we discuss general invariance theorems for
the sharp constants in the
Markov-Bernstein-Nikolskii type inequalities.

It is easy to see that the sharp constant $C_n^*$ in
inequality \eqref{E1.1} for $p\in [1,\iy]$ does not
change if one restricts \eqref{E1.1} to even polynomials $P$
(i.e., $P(x)=P(-x),\,\vert x\vert\le 1).$ The proof of this fact involves
first, the symmetrization operator $\Tilde{P}(x):=(P(x)+P(-x))/2$ of $P$
and second, the invariance of the operator $\Delta^N$ and the
class $\PP_{n,m}$ under the group $G_m=\{-e^*,e^*\}$, where $e^*$ is
the identity transformation on the unit ball $\Om$. Indeed,
$\vert \Delta^N\Tilde{P}(0)\vert=\vert \Delta^NP(0)\vert$ and
\ba
\left(\int_{\vert x\vert\le 1}\vert \Tilde{P}(x)\vert^pdx\right)^{1/p}
\le \left(\int_{\vert x\vert\le 1}\vert P(x)\vert^pdx\right)^{1/p},
\qquad p\in [1,\iy].
\ea
Below we extend this example to compact topological groups $G_m$
of continuous transformations on more general sets $\Om$ and to
more general classes of elements. The symmetrization will be provided
 by the Haar integral. More examples are discussed in Section \ref{S3}.

Let $\Om \subseteq \R^m$ and let $F(\Om)$ be a Banach space of functions
 $f:\Om\to \CC^1$ with the norm $\|\cdot\|_{F(\Om)}$.
 Next, let $\B \ne  \{0\}$ be a closed subspace of $F(\Om)$ and let
 $D:\B\to F(\Om)$ be a bounded linear operator.
 Given $a\in\Om$, we define the sharp constant in
 the generalized  Markov-Bernstein-Nikolskii type inequality by
 \beq\label{E2.0}
\C_a= \C_a\left( \B, F(\Om), D\right):=\sup_{f\in \B\setminus \{0\}}
 \frac{\left\vert D(f)(a)\right\vert}
 {\|f\|_{F(\Om)}}.
 \eeq
 If $D$ is the imbedding operator $I:\B\to F(\Om)$, then $\C_a$ is
 the sharp constant in
 the generalized  Nikolskii-type inequality.

 Further, let $G_m=G_m(a)$ be a compact topological group of
 continuous transformations $s:\Om\to\Om$ with a fixed point $a\in\Om$
 (i. e., $sa=a,\,s\in G_m$) and let $F(\Om)^{G_m}$ denote a subspace of
 $F(\Om)$ of all functions $f$ which are invariant under the group $G_m$, i.e.,
 $f(s\cdot) =f(\cdot),\,s\in G_m$.

 Let $B$ be a closed subspace of $F(\Om)$. In this section we discuss sufficient
  conditions for the sharp constant to be invariant under $G_m$, i.e.,
 \beq\label{E2.1}
 \C_a\left( B, F(\Om), D\right)=\C_a\left( B\cap F(\Om)^{G_m}, F(\Om), D\right).
 \eeq

 We assume that $B,\,G_m,\,D,$ and $F(\Om)$ satisfy the following conditions.
 \begin{itemize}
 \item[(C1)] The norm $\|\cdot\|_{F(\Om)}$ is invariant under $G_m$,
 i.e., for every $f\in F(\Om)$ and each $s\in G_m,\,
 \|f(s\cdot)\|_{F(\Om)}=\|f\|_{F(\Om)}$.
 \item[(C2)] The operator $D$ is invariant under $G_m$,
 i.e., for every $f\in B$ and each $s\in G_m,\,D(f_s)(\cdot)=D(f)(s\cdot)$,
 where $f_s(\cdot):=f(s\cdot)$.
 \item[(C3)] The subspace $B$ is invariant under $G_m$,
 i.e., for every $f\in B$ and each $s\in G_m,\,f(s\cdot)\in B$.
 \item[(C4)] For every $f\in B$ the mapping $f(s\cdot):G_m\to B$ is a continuous
 function in $s\in G_m$, i.e., for any $s\in G_m,\,
 \lim_{\tau\to s}\|f(s\cdot)-f(\tau\cdot)\|_{F(\Om)}=0$.
 \item[(C5)] $B\cap F(\Om)^{G_m} \ne \{0\}$.
 \end{itemize}
 The following general invariance theorem holds true.

 \begin{theorem}\label{T2.1}
 If conditions (C1) through (C5) are satisfied, then \eqref{E2.1} is valid.
 \end{theorem}
 \noindent
 \proof
 It suffices to prove the inequality
 \beq\label{E2.2}
 \frac{\left\vert D(f)(a)\right\vert}{\|f\|_{F(\Om)}}
 \le \sup_{f^*\in \left(B\cap F(\Om)^{G_m}\right)\setminus \{0\}}
 \frac{\left\vert D(f^*)(a)\right\vert}
 {\|f^*\|_{F(\Om)}}
 \eeq
 for every $f\in B\setminus\{0\}$. Note that due to condition (C5)
 the right-hand sides of \eqref{E2.1} and \eqref{E2.2} are well-defined.

 Since $G_m$ is a compact topological group, there exists the Haar
 measure $\mu(s)$ on $G_m$ with $\mu(G_m)=1$
 (see, e.g., \cite[Theorem 5.14]{R1973}). Next, for every $f\in B$
 the function $f(s\cdot):G_m\to B$ is continuous on $G_m$ by conditions
 (C3) and (C4); therefore, its image $H:=\{f(s\cdot): s\in G_m\}$
 is compact in $F(\Om)$. Since $F(\Om)$ is a Banach space, the closure of
 the convex hull of $H$ (denote it by $\Bar{H}_{co}$) is compact as well
 (see, e.g., \cite[Theorem 3.25 (a)]{R1973}).
 Then the Haar integral
 \beq\label{E2.3}
 f^*(\cdot):=\int_{G_m}f(s\cdot)d\mu(s)
 \eeq
 exists and $f^*\in \Bar{H}_{co}$ (see, e.g. \cite[Theorem 3.27]{R1973}).
 Note that examples of the Haar measures and the Haar integral can be found
 in \cite [Sect. 15.17]{HR1963} and \cite[Sect. 12.1]{S1976}.
 Moreover, since $B$ is a closed subset of $F(\Om)$, we conclude that
 $\Bar{H}_{co} \subseteq B$, so $f^*\in B$.
 Next, for every $t\in G_m$,
 \beq\label{E2.4}
 f^*(t\cdot)=\int_{G_m}f(st\cdot)d\mu(s)=\int_{G_m}f(s\cdot)d\mu(s)=f^*(\cdot),
 \eeq
 where the second equality in \eqref{E2.4} follows from the invariance of the
 Haar measure. Therefore, $f^*\in B\cap F(\Om)^{G_m}$. Using now the generalized
  Minkowski inequality (see, e.g., \cite[Lemma 3.2.15]{DS1988}) and condition (C1),
  we obtain
  \beq\label{E2.5}
  \|f^*\|_{F(\Om)}\le \int_{G_m}\|f(s\cdot)\|_{F(\Om)}d\mu(s) =\|f\|_{F(\Om)}.
  \eeq
  Further, by \cite[Theorem 3.2.19(c)]{DS1988} and conditions (C3) and (C2), we have
  \beq\label{E2.6}
  D(f^*)(x)=\int_{G_m}D(f_s)(x)d\mu(s)=\int_{G_m}D(f)(s x)d\mu(s),\qquad x\in\Om,
  \eeq
  and since $sa=a,\,s\in G_m$, we obtain from \eqref{E2.6}
  \beq\label{E2.7}
  D(f^*)(a)=D(f)(a).
  \eeq
  If $f\ne 0$ and $D(f)(a)=0$, then \eqref{E2.2} holds trivially true by (C5).
  If $D(f)(a)=D(f^*)(a)\ne 0$, then $f\ne 0,\,f^*\ne 0$, and it follows from
  \eqref{E2.5} and \eqref{E2.7} that
  \ba
  \frac{\left\vert D(f)(a)\right\vert}{\|f\|_{F(\Om)}}
  \le \frac{\left\vert D(f^*)(a)\right\vert}{\|f^*\|_{F(\Om)}}.
  \ea
  Hence \eqref{E2.2} holds true in this case as well. Thus \eqref{E2.1}
  is established. \hfill $\Box$\vspace{.12in}\\
  The proof of Theorem \ref{T2.1} is based on the existence of the Haar integral
  $f^*$ which belongs to $B$. These both facts follow from strong condition (C4).
  However, the existence of $f^*(x)$ for each $x\in\Om$ follows from the
  following weaker condition.
  \begin{itemize}
  \item[(C4$^\prime$)] For every $f\in B$ and each fixed $x\in\Om$, the linear
  functional $f(sx):G_m\to \CC^1$ is continuous in $s\in G_m$.
  \end{itemize}
  It is obvious that (C4) implies (C4$^\prime$), but the converse statement is not
   valid in general (see \cite[Example 2.1]{G2008}). If we introduce a new
   condition
  \begin{itemize}
  \item[(C6)] For every $f\in B$, the Haar integral $f^*$ defined in \eqref{E2.3}
   belongs to $B$,
   \end{itemize}
   then we arrive at the following version of Theorem \ref{T2.1}.

   \begin{theorem}\label{T2.2}
   If conditions (C1), (C2), (C3), (C4$^\prime$), (C5), and (C6) are satisfied,
   then \eqref{E2.1} is valid.
   \end{theorem}
   \proof
   The existence of $f^*$ follows from (C4$^\prime$)
   (see, e.g., \cite[Theorem 5.14]{R1973}) and, in addition,
   $f^*\in B$ by (C6). The rest of the proof of \eqref{E2.1} is
   similar to that of Theorem \ref{T2.1}.\hfill $\Box$

   Note that invariant means like the Haar integral exist
   in more general situation, more precisely, on amenable semigroups
    (see \cite[Sect. 6]{A1996} for details).

   \section{Special Invariance Theorems}\label{S3}
 \noindent
\setcounter{equation}{0}
Here, we discuss special cases of invariance theorems presented in
Section \ref{S2}.\vspace{.12in}\\
   \textbf{Special Cases and Preliminaries.}
   In all our examples of applications of Theorems \ref{T2.1} and \ref{T2.2}
   we  use special sets $\Om$, spaces $F(\Om)$, subspaces $B$,
   groups $G_m$, and linear operators $D$ described below.
   In addition, we  discuss here their certain properties.

   Let $\Om$ be one of the following sets: $\BB^m,\,m\ge 1;\,
   S^{m-1},\,m\ge 2;$ and $\R^m,\,m\ge 1$; and let
   $F(\Om)=L_p(\Om),\,1\le p\le\iy$, be the space of all measurable functions
   $f:\Om\to\CC^1$ with the finite norm
   \ba
 \|f\|_{L_p(\Om)}:=\left\{\begin{array}{ll}
 \left(\int_\Om\vert f(x)\vert^p dx\right)^{1/p}, & 1\le p<\iy,\\
 \mbox{ess} \sup_{x\in \Om} \vert f(x)\vert, &p=\iy,
 \end{array}\right.
 \ea
 if $\Om=\BB^m$ or $\Om=\R^m$ and
 \ba
 \|f\|_{L_p(S^{m-1})}:=\left\{\begin{array}{ll}
 \left(\int_{S^{m-1}}\vert f(x)\vert^p dS(x)\right)^{1/p}, & 1\le p<\iy,\\
 \mbox{ess} \sup_{x\in S^{m-1}} \vert f(x)\vert, &p=\iy,
 \end{array}\right.
 \ea
 where $S(\cdot)$ is the spherical surface Lebesgue measure
 on $S^{m-1}$.
 In addition, we also need the weighted space
 $L_{p,\mu(t)}(\Omega_1),\,0<p\le\iy,$
 of all univariate measurable functions
   $f:\Omega_1\to\CC^1$ with the finite quasinorm
   \ba
 \|f\|_{L_{p,\mu(t)}(\Omega_1)}:=\left\{\begin{array}{ll}
 \left(\int_{\Omega_1}\vert f(t)\vert^p\mu(t) dt\right)^{1/p},
  & 0< p<\iy,\\
 \mbox{ess} \sup_{t\in \Omega_1} \vert f(t)\vert, &p=\iy.
 \end{array}\right.
 \ea
 Here, $\Omega_1$ is a measurable subset of  $\R^1$  and
 $\mu:\Omega_1\to [0,\iy)$ is a locally integrable weight.
 This quasinorm allows the following "triangle" inequality
 \beq\label{E2.7a}
 \left\|f+g\right\|^{\tilde{p}}_{L_{p,\mu(t)}(\Omega_1)}
 \le \left\|f\right\|^{\tilde{p}}_{L_{p,\mu(t)}(\Omega_1)}
 +\left\|g\right\|^{\tilde{p}}_{L_{p,\mu(t)}(\Omega_1)},
 \eeq
 where $\tilde{p}:=\min\{1,p\}$ for $p\in(0,\iy]$.
 In this section, $\mu(t)$
 is either $\vert t\vert^{m-1},\,m\ge 1$, or $(1-t^2)^{(m-3)/2},\,m\ge 2$.
 In Sections \ref{S4}, \ref{S5}, and \ref{S6}, we use more general weights.

 In the capacity of $B$ we discuss either the set
  $\PP_{n,m}\left\vert_{\Om}\right.$ of the restrictions
  $P\left\vert_{\Om}\right.$ to $\Om$ of
  polynomials $P(x)=\sum_{\vert\al\vert\le n}c_\al x^\al\,:\R^m\to\CC^1$
   in $m$ variables of degree at most $n$
  (if $\Om$ is identified, we often write  $\PP_{n,m}$ instead of
    $\PP_{n,m}\left\vert_{\Om}\right.$) or the set
  $E_{\sa,m}\cap L_p(\R^m)$ of the restrictions to $\R^m$ of entire functions
  in $m$ variables of spherical type $\sa>0$ that belong to $L_p(\R^m)$.

  We recall that an entire function $f:\CC^m\to\CC^1$ has spherical type
  $\sa>0$ if for any $\vep>0$ there exists a constant $C_0(\vep,f)$ such that
  \beq\label{E2.8}
  \vert f(z)\vert\le C_0(\vep,f)\exp(\sa(1+\vep)\vert z\vert),\qquad z\in\CC^m,
  \eeq
  (see \cite[Sect. 3.2.6]{N1969} and \cite[Definition 5.1]{DP2010}).
  We  often identify $f$ with its restriction to $\R^m$. Note that
  $E_{\sa,m}\cap L_p(\R^m)$ is a closed subspace of $L_p(\R^m),\,p\in(0,\iy]$
  (cf. \cite[Theorem 3.5]{N1969}).

  We need the following compactness theorem for functions from
  $E_{\sa,m}\cap L_p(\R^m)$.

  \begin{proposition}\label{P2.3}
  For any sequence $\{f_n\}_{n=1}^\iy,\,
f_n\in E_{\sa,m}\cap L_p(\R^m),\,n\in\N,\,p\in[1,\iy]$,
with $\sup_{n\in\N}\| f_n\|_{L_\iy(\R^m)}= C$, there exist a subsequence
$\{f_{n_s}\}_{s=1}^\iy$ and a function $f_0\in E_{\sa,m}\cap L_p(\R^m)$
such that
\beq\label{E2.9}
\lim_{s\to\iy} f_{n_s}=f_0
\eeq
uniformly on any compact set in $\CC^m$.
\end{proposition}
\proof
It follows from \eqref{E2.8} that if $f\in E_{\sa,m}$, then
\beq\label{E2.9a}
\vert f(z)\vert\le C_0(\vep,f)\exp\left(\sa(1+\vep)
\sum_{j=1}^m\vert z_i\vert\right),\qquad z\in\CC^m.
\eeq
Therefore, $f_n,\,n\in\N,$ has exponential type $\sa$ by the definition in
\cite[Sect. 3.1]{N1969}.
Since $f_n\in E_{\sa,m}\cap L_p(\R^m)$,
 by Nikolskii's compactness theorem \cite[Theorem 3.3.6]{N1969},
 there exists a subsequence
$\{f_{n_s}\}_{s=1}^\iy$ and an entire function $f_0\in L_\iy(\R^m)$
such that \eqref{E2.9} holds true uniformly on any compact set in $\CC^m$.
In addition, $f_0\in E_{\sa,m}.$ Indeed, since
\beq\label{E2.10}
\vert f_n(x+iy)\vert\le C\exp(\sa\vert y\vert),
\qquad n\in\N,\quad x\in\R^m,\quad y\in\R^m,
\eeq
(see \cite[Eq. (4.13)]{NW1978}), we obtain by \eqref{E2.10}
\ba
\vert f_0(x+iy)\vert=\lim_{s\to\iy}\vert f_{n_s}(x+iy)\vert
\le C \exp(\sa\vert y\vert)\le C\exp(\sa\vert z\vert).
\ea
Thus $f_0\in E_{\sa,m}\cap L_p(\R^m)$. \hfill $\Box$\vspace{.12in}\\

In addition to $\PP_{n,m}=\PP_{n,m}\left\vert_{\Om}\right.$  and $E_{\sa,m}$,
we  also need univariate sets
$\PP_{n,1,e}=\PP_{n,1,e}\left\vert_{[-1,1]}\right.$ and $E_{\sa,1,e}$ of all
even polynomials and even entire functions from
$\PP_{n,1}=\PP_{n,1}\left\vert_{[-1,1]}\right.$ and $E_{\sa,1}$, respectively.

Throughout the paper we use the following groups $G_m=G_m(a)$.
Let $G_m(0)=O(m)$ be the group of all  proper and improper rotations
(about the origin) of $\R^m$.
We identify $O(m)$ with the group of all $m \times m$ orthogonal matrices
which is isomorphic to $O(m)$ since $s\in O(m)$ if and only if
$sx=A(s)\,x^T$, where $A(s)$ is an $m \times m$ orthogonal matrix with
$\vert \det A(s)\vert=1$ and $x^T$ is a column vector.
Let $G_m(a)=O(m,a)$ be a subgroup of $O(m)$ of all proper and
 improper rotations (or
$m \times m$ orthogonal matrices) $s$, satisfying the condition $sa=a$,
where $a\ne 0$ is a fixed vector from $\R^m$. For example, if
$a= (\cos \g, \sin \g)$, then $O(2,a)=\{I, A_\g\}$, where $I$ is the
 $2 \times 2$ identity matrix and
 \ba
 A_\g=\begin{bmatrix}
 \cos 2\g &\sin 2\g\\
 \sin 2\g & -\cos 2\g
 \end{bmatrix},
 \ea
 which is the product of two transformations: reflection about the $x$-axis
 and rotation by the angle $2\g$.
 Note that $O(m)$ will be used in Examples \ref{Ex2.10} and \ref{Ex2.12},
 while $O(m,a)$ will be used in Example \ref{Ex2.14}.

Throughout the paper we use the polyharmonic operator $D=\Delta^N$,
 where $N\in\Z^1_+$ and
 \ba
 \Delta=\Delta_x:=\sum_{j=1}^m\frac{\partial^2}{\partial x_j^2}
 \ea
 is the Laplace operator on $\Om$.
 Note that the restriction $\boldsymbol{\de}_x$
  of $\Delta_x$ on $S^{m-1}$ is defined by
 $\boldsymbol{\de}_x(f)(x)=\Delta_x(f)\left\vert_{S^{m-1}}\right.(x)
 :=\Delta_x(f)(x/\vert x\vert)$.
 In case of $N=0,\,\Delta^N$ is the imbedding operator
 $I:\PP_{n,m}\left\vert_{\Om}\right.\to L_p(\Om)$
 or $I:E_{\sa,m}\cap L_p(\R^m)\to L_p(\R^m)$.
  We need certain properties of $\Delta$.
 \begin{proposition}\label{P2.4}
 Let $\R^m,\,m\ge 2$, be equipped with the spherical coordinates
 $x=(r,\,\theta_1,\ldots,\theta_{m-1}),\,r\in[0,\iy),\,\theta_j\in[0,\pi],\,
 1\le j\le m-2,\, \theta_{m-1}\in[0,2\pi).$ Then the following properties
 of $\Delta$ hold true.\\
 (a) In spherical coordinates,
 \ba
 \Delta_x=\frac{\partial^2}{\partial r^2}+\frac{m-1}{r}\frac{\partial}{\partial r}
 +\frac{1}{r^2}\boldsymbol{\de},
 \ea
 where $\boldsymbol{\de}$ is the spherical Laplacian given by
 \ba
 \boldsymbol{\de}:=\sum_{j=1}^{m-1}\frac{1}{q_j\sin^{m-j-1}\theta_j}
 \frac{\partial}{\partial\theta_j}\left(\sin^{m-j-1}\theta_j
 \frac{\partial}{\partial\theta_j}\right)
 \ea
 and $q_1:=1; q_j:=\prod_{s=1}^{j-1}\sin^2\theta_s,\,2\le j\le m-1$.
 In particular, $\boldsymbol{\de}$ is the representation of $\boldsymbol{\de}_x$
 in spherical coordinates.\\
 (b) For a fixed $a\in S^{m-1}, \,l\in\R^1,$ and $k\in\R^1,$
 \ba
 \Delta_x\left(r^l(x,a)^k\right)=\left(\frac{k(k-1)}{(x,a)^2}+
 \frac{l(l+m+2k-2)}{r^2}\right)r^l(x,a)^k,\qquad r=\vert x\vert.
 \ea
 (c) If $f(x)$ is a radial function $\vphi(\vert x\vert)$,
 where $\vphi:\R^1\to\CC^1$ is an even twice continuously differentiable
  function on $\R^1$, then
  $\Delta_x(f)(x)=Be_{m/2-1}\left(\vphi\right)(\vert x\vert),\,m\ge 1$.
  Here,
  \beq\label{E2.11}
  Be_{\nu}(\vphi)(r)
  :=\vphi^{\prime\prime}(r)+\frac{2\nu+1}{r}
  \vphi^\prime(r),\qquad \nu\ge -1/2,\quad r\in\R^1,
  \eeq
  is the Bessel operator
  and $Be_\nu(\vphi)(0):=\lim_{r\to 0}Be_\nu(\vphi)(r)$.\\
  (d) If $f(x)=\vphi((x,a))$, where
  $\vphi\in\PP_{n,1},\,x\in S^{m-1}$,
  and $a\in S^{m-1}$ is a fixed point, then
  $\boldsymbol{\de}_x(f)(x)=Ge_{m/2-1}\left(\vphi\right)((x,a)),\,m\ge 2$.
  Here,
  \beq\label{E2.12}
  Ge_{\la}\left(\vphi\right)(t)
  :=(1-t^2)\vphi^{\prime\prime}(t)-(2\la+1)\, t\,
  \vphi^\prime(t),\qquad \la\ge -1/2,\quad t\in[-1,1],
  \eeq
  is the Gegenbauer operator.\\
  (e) The operators $\Delta_x$ and $\boldsymbol{\de}_x$
   are invariant under orthogonal transformations,
   i.e., for $s\in O(m)$ and $y=sx,\,\Delta_x=\Delta_y$
   and $\boldsymbol{\de}_x=\boldsymbol{\de}_y$.
 \end{proposition}
 \proof
 Statements (a), (b), and (e) are well-known and can be found in
 \cite[Sect. 11.1.1]{Erd1953}, while (c) follows immediately from (a).
 It suffices to prove (d) for $\vphi(t)=t^k,\,k=0,\,1,\,\ldots,n.$
 Using statement (b) for $l=-k$ and $r=1$, we obtain
 \ba
 \boldsymbol{\de}_x((x,a)^k)=\Delta_x(r^{-k}(x,a)^k)
 =k(k-1)(x,a)^{k-2}-k(k+m-2)(x,a)^k
 =Ge_{m/2-1}(\vphi)(t)\left\vert_{t=(x,a)}\right..
 \ea
 This completes the proof of the proposition. \hfill $\Box$
 \begin{remark}\label{R2.5}
 We call \eqref{E2.11} the Bessel operator because the functions
 $t^{-\nu}J_{\nu}(\sqrt{c} t)$ are eigenfunctions of $Be_{\nu}$
 for $c\ge 0$, see \cite[Sect. 4.31]{W1944}.
 We call \eqref{E2.12} the Gegenbauer operator because the
 Gegenbauer polynomials $C_k^{\la}$
  are eigenfunctions of $Ge_{\la}$, see \cite[Sect. 10.9]{Erd1953}.
 \end{remark}
 \noindent
 \textbf{Conditions in Special Cases.}
 We first discuss  conditions (C1), (C2), and (C3) in special cases.

 \begin{proposition}\label{P2.6}
 Let $G_m=G_m(a)$ be a subgroup of $O(m,a)$. Then the following
 statements hold true.\\
 (a) If $\Om$ is one of the sets $\BB^m,\,m\ge 1;\,S^{m-1},m\ge 2;
 \,\R^m,\,m\ge 1$, and
 $F(\Om)=L_p(\Om),\,1\le p\le \iy$, then (C1) is satisfied.\\
 (b) If $D$ is the polyharmonic operator $\Delta^N,\,N=0,\,1,\dots,$
 then (C2) is satisfied.\\
 (c) If $B=\PP_{n,m}\left\vert_{\Om}\right.$, where either
 $\Om=\BB^m,\,m\ge 1$, or
 $\Om=S^{m-1},\,m\ge 2$, then (C3) is satisfied.\\
 (d)  If $B=E_{1,m}\cap L_p(\R^m)$,  then (C3) is satisfied.
 \end{proposition}
 \proof
 Statement (a) is obviously satisfied, while (b) follows from Proposition
 \ref{P2.4} (e). Let $P\in \PP_{n,m}\left\vert_{\Om}\right.$ and
  $x\in\Om$.
 Since $sx\in\Om$, we have
 \beq\label{E2.13}
 P\left\vert_{\Om}\right.(sx)=P(sx)=P(s\cdot)\left\vert_{\Om}\right.(x).
 \eeq
 Next, $s$ is a linear transformation, so $P(s\cdot)\in\PP_{n,m}
 \left\vert_{\R^m}\right.$ and
 $P\left\vert_{\Om}\right.(s\cdot)\in \PP_{n,m}\left\vert_{\Om}\right.$
 by \eqref{E2.13}.

 Further, let $f\in E_{1,m}\cap L_p(\R^m)$. Extending $f$ and
  $f(s\cdot)$ to $\CC^m$,
  we see that $f(s\cdot)$ is an entire function since $s\in O(m)$ is a
  linear transformation. Moreover, for
  $z\in\CC^m,\,\vert sz\vert=\vert z\vert$.
  Therefore, $f(s\cdot)\in E_{1,m}$ by \eqref{E2.8}.
  In addition, $f(s\cdot)\in L_p(\R^m)$ by statement (a).
  Thus statements (c) and (d) are established.\hfill $\Box$\vspace{.12in}\\

  In the following two propositions we discuss the validity of
  conditions (C4) and (C4$^\prime$) in special cases.
  \begin{proposition}\label{P2.7}
  Let $\Om$ be a closed subset of $\R^m$ and $F(\Om)=L_p(\Om),
  \,1\le p\le\iy$. Next, let $B$ be a subspace of $L_p(\Om)$ of
  continuous functions on $\Om$, and let $G_m^*$ be a compact group of
  linear transformations of the form $sx=A(s)\, x^T:\Om\to\Om$,
  where $A(s)$ is an $m \times m$ matrix. Then the following
   two statements hold true.\\
   (a) Condition (C4$^\prime$) is satisfied for $1\le p\le\iy$.\\
   (b) Condition (C4) is satisfied for $1\le p\le\iy$ and a
   compact set $\Om$.
   \end{proposition}
   \proof
   Statement (a) follows from the uniform continuity of $f(sx)$
   in $s\in G_m^*$ for each $x\in\Om$ if we take into account
   the elementary inequality
   \ba
   \vert sx-s_1x\vert\le \|A(s)-A(s_1)\|_F \vert x\vert,
   \qquad s\in G_m^*,\quad s_1\in G_m^*,\quad x\in \Om.
   \ea
   Similarly, statement (b) follows from the uniform continuity
    of $f(s\cdot)$
   in $s\in G_m^*$, the estimate
   \ba
   \max_{x\in\Om}\vert sx-s_1x\vert\le \|A(s)-A(s_1)\|_F
   \max_{x\in\Om}\vert x\vert,
   \qquad s\in G_m^*,\quad s_1\in G_m^*,
   \ea
   and the continuous imbedding of $L_\iy(\Om)$ into $L_p(\Om)$ for
   a compact set $\Om$.\hfill $\Box$
   \begin{proposition}\label{P2.8}
   Let $B$ be a  subspace of $L_p(\R^m),\,1\le p<\iy,$ of
  continuous functions on $\R^m$, and let $G_m^*$ be a subgroup of $O(m)$.
  Then condition (C4) is satisfied.
  \end{proposition}
  \proof
  Let $f\in B$ and let $s$ and $s_1$ be two proper or improper rotations
  and let $A(s)$ and $A(s_1)$ be the corresponding orthogonal matrices.
  Given $\vep>0$ there exists $M=M(\vep,f)>0$ such that
  \beq\label{E2.14}
  \left(\int_{\vert x\vert>M}\vert f(sx)\vert^pdx\right)^{1/p}
  = \left(\int_{\vert x\vert>M}\vert f(s_1x)\vert^pdx\right)^{1/p}
  = \left(\int_{\vert x\vert>M}\vert f(x)\vert^pdx\right)^{1/p}
  <\vep/3.
  \eeq
  Next, by statement (b) of Proposition \ref{P2.7}, there exists
  $\de(\vep)>0$ such that for $\|s-s_1\|:=\|A(s)-A(s_1)\|_F<\de(\vep)$,
   \beq\label{E2.15}
   \|f(s\cdot)-f(s_1\cdot)\|_{L_p(\BB^m(M))}<\vep/3.
   \eeq
   Combining \eqref{E2.14} with \eqref{E2.15}, we obtain
   \ba
    \|f(s\cdot)-f(s_1\cdot)\|_{L_p(\R^m)}<2\vep/3+\vep/3=\vep.
    \ea
    Then condition (C4) is satisfied.\hfill $\Box$
    \begin{remark}\label{R2.9}
    Note that condition (C4) is not always satisfied if
    $B$ is a  subspace of $L_\iy(\R^m)$ of
  continuous functions on $\R^m$,
    see \cite [Example 2.1]{G2008}.
    \end{remark}
    \noindent
    \textbf{Examples.}
    Here, we discuss typical examples of applications of Theorems
     \ref{T2.1} and \ref{T2.2}.
     In addition to $D=\Delta^N$, we also use Bessel and Gegenbauer
     operators $D=(Be_\nu)^N$ and $D=(Ge_\la)^N$ defined
     by \eqref{E2.11} and \eqref{E2.12}.
     In case of $N=0,\, (Be_\nu)^N$ and $(Ge_\la)^N$ are
     the corresponding imbedding operators.
     \begin{example}\label{Ex2.10}
     $\Om=\BB^m,\,F(\Om)=L_p(\BB^m),\,m\ge 1,\,1\le p\le\iy;\,
     G_m=G_m(0)=O(m);
     \linebreak L_p(\BB^m)^{O(m)}$ is the set of all radial
     functions from $L_p(\BB^m);\,B=\PP_{n,m};\,
     D=\Delta^N=\Delta^N\left\vert_{\BB^m}\right.,\linebreak
     N\in\Z^1_+$.

     Conditions (C1), (C2), and (C3) are satisfied by Proposition
     \ref{P2.6} and (C4) is satisfied by Proposition
     \ref{P2.7} (b).
     In addition, for $1\le p\le\iy$,
     \beq\label{E2.16}
     \PP_{n,m}\cap L_p(\BB^m)^{O(m)}
     =\left\{P(x)=Q(\vert x\vert):
     Q\in \PP_{2\lfloor n/2 \rfloor, 1,e}\right\}
     \eeq
     (see \cite[Lemma 4.2.11]{SW1971} and \cite[Proposition 4.1]{G2008}).
     In particular, condition (C5) is satisfied by \eqref{E2.16}.

     Using \eqref{E2.16} and Proposition
     \ref{P2.4} (c), we obtain from Theorem \ref{T2.1}
     \ba
     \sup_{P\in\PP_{n,m}\setminus\{0\}}
     \frac{\left\vert\Delta^N (P)(0)\right\vert}
     {\left\|P\right\|_{L_p(\BB^m)}}
     =\left(\frac{2}{\left\vert S^{m-1}\right\vert_{m-1}}\right)^{1/p}
     \sup_{Q\in\PP_{2\lfloor n/2 \rfloor, 1,e}
     \setminus\{0\}}
     \frac{\left\vert \left(Be_{m/2-1}\right)^N (Q)(0)\right\vert}
     {\left\|Q\right\|_{L_{p,\vert t\vert^{m-1}}([-1,1])}}.
     \ea
     Thus we arrive at the following corollary.
     \begin{corollary}\label{C2.11}
     For $m\in\N,\,n\in\N,\,N\in\Z^1_+$, and $p\in[1,\iy]$,
     \bna\label{E2.17}
    && \C_0\left( \PP_{n,m},\,L_p(\BB^m),\,
     \Delta^N\right)\nonumber\\
     &&=\left(\frac{2}{\left\vert S^{m-1}\right\vert_{m-1}}\right)^{1/p}
     \C_0\left(\PP_{2\lfloor n/2 \rfloor, 1,e}
    ,\,L_{p,\vert t\vert^{m-1}}([-1,1]),\,
     \left(Be_{m/2-1}\right)^N\right).
     \ena
     \end{corollary}
     \end{example}

     \begin{example}\label{Ex2.12}
     $\Om=\R^m,\,F(\Om)=L_p(\R^m),\,m\ge 1,\,1\le p\le\iy;\,
     G_m=G_m(0)=O(m);\linebreak L_p(\R^m)^{O(m)}$ is the set of all radial
     functions from $L_p(\R^m);\,B=E_{1,m}\cap L_p(\R^m);\,
     D=\Delta^N,\,N\in\Z^1_+$.

     Conditions (C1), (C2), and (C3) are satisfied by Proposition
     \ref{P2.6} and  condition
     (C4) for $p\in[1,\iy)$ is satisfied by Proposition
     \ref{P2.8}.

     In case of $p=\iy$, condition (C4$^\prime$) is satisfied
     by Proposition
     \ref{P2.7} (a).
     Next, we show that condition (C6) is satisfied for $p=\iy$
     and condition (C5) is satisfied for $1\le p\le\iy$.

     We first state the following simple fact
     (cf. \cite[Lemma 3.6.1]{N1969}) that is used in this example:
     if $h(t)\in E_{\sa,1}$ is an even function, then
     $\vphi(x):=h(\vert x\vert)\in E_{\sa,m}$. Indeed, it is an entire function
     with the extension $h(\vert z\vert)$ to $\CC^m$. Moreover, since
     $\vert h(\vert z\vert)\vert
     \le C_0(\vep,h)\exp(\sa(1+\vep)\vert z\vert)$ by\eqref{E2.8},
     we conclude that $\vphi\in E_{\sa,m}$.

     In particular, the function
     \ba
     h_n(x):=\left(\frac{\sin\left(\frac{\vert x\vert}
     {n(m+1)}\right)}{\vert x\vert/[n(m+1)]}\right)^{m+1}
     \ea
     belongs to $E_{1/n,m},\,n\in\N.$

     Now we show that condition (C6) is satisfied, i.e.,
     for $f\in E_{1,m}\cap L_\iy(\R^m)$ the Haar integral
     \beq\label{E2.17a}
     f^*(x)=\int_{s\in O(m)}f(sx)d\mu(s)
     \eeq
     belongs to $E_{1,m}\cap L_\iy(\R^m)$. Let us set
     $f_n(x):=f(x)h_n(x),\,n\in \N$.
     Then \linebreak $f_n\in E_{1+1/n,m}\cap L_1(\R^m),\,n\in\N$, and by
     the elementary inequality $1-\sin t/t\le t^2/6,\,t\in R^1$,
     we obtain
     \beq\label{E2.18}
     \vert f(x)-f_n(x)\vert\le \frac{\vert x\vert^2}
     {6n^2(m+1)}\|f\|_{L_\iy(\R^m)},\qquad n\in\N,\quad x\in\R^m.
     \eeq
     Next, for the corresponding Haar integral we have
     \beq\label{E2.19}
     f^*_n(x)=\int_{s\in O(m)}f_n(sx)d\mu(s)\in E_{1+1/n,m}\cap L_1(\R^m).
     \eeq
     We prove \eqref{E2.18} similarly to the fact that $f^*\in B$
     in the proof of Theorem \ref{T2.1}. The mapping $f_n(s\cdot)
     :O(m)\to E_{1+1/n,m}\cap L_1(\R^m)$
     is continuous on $O(m)$ by condition (C3) (i.e.,
     $f_n\in E_{1+1/n,m}\cap L_1(\R^m)$; see Proposition \ref{P2.6} (d))
     and by Condition (C4) (see Proposition  \ref{P2.7} (b)).
     Then its image $H:=\{f_n(s\cdot):s\in O(m)\}$ is compact in $L_1(\R^m)$
     and the closure of the convex hull $\Bar{H}_{co}$ of $H$ is compact as well
     (see, e.g., \cite[Theorem 3.25 (a)]{R1973}).
     Then $f_n^*$ exists and $f_n^*\in \Bar{H}_{co}$; since
     $\Bar{H}_{co}\subseteq E_{1+1/n,m}\cap L_1(\R^m)$, we conclude that
     $f_n^*\in  E_{1+1/n,m}\cap L_1(\R^m)$.

     In addition, it follows from \eqref{E2.17a}, \eqref{E2.18}, and \eqref{E2.19}
     \bna\label{E2.20}
    \vert f^*(x)-f^*_n(x)\vert
    &\le&
     \int_{s\in O(m)}\vert f(sx)- f_n(sx)\vert d\mu(x)\nonumber\\
     &\le& \frac{\vert sx\vert^2}
     {6n^2(m+1)}\|f\|_{L_\iy(\R^m)}
     = \frac{\vert x\vert^2}
     {6n^2(m+1)}\|f\|_{L_\iy(\R^m)}.
     \ena
     Next, \eqref{E2.20} shows that
     \beq\label{E2.21}
     \lim_{n\to\iy}f_n^*=f^*
     \eeq
     uniformly on any compact set of $\R^m$.
     Further,
     \ba
     \sup_{n\in\N}\left\|f_n^*\right\|_{L_\iy(\R^m)}\le \|f\|_{L_\iy(\R^m)}.
     \ea
     Then by Proposition \ref{P2.3} and \eqref{E2.19}, there exist a subsequence
     $\{f^*_{n_s}\}_{s=1}^\iy$ and a function
     $f_0^*\in E_{1,m}\cap L_\iy(\R^m)$ such that
     \beq\label{E2.22}
     \lim_{s\to\iy}f_{n_s}^*(x)=f_0^*(x)
     \eeq
     uniformly on any compact set of $\R^m$. Comparing \eqref{E2.21}
     and \eqref{E2.22}, we conclude that $f^*=f_0^*$, so
     $f^*\in E_{1,m}\cap L_\iy(\R^m)$. Thus condition (C6) is satisfied
     in the case $p=\iy$.

     In addition, for $1\le p\le\iy$,
     \beq\label{E2.23}
     E_{1,m}\cap L_p(\R^m)^{O(m)}
     =\left\{f(x)=g(\vert x\vert):
     g\in E_{1,1,e}\cap L_{p,\vert t\vert^{m-1}}(\R^1)\right\}
     \eeq
     (see \cite[Proposition 6.1]{G2008}).
     For example, the function
     \ba
     h_1(x):=\left(\frac{\sin\left(\frac{\vert x\vert}
     {m+1}\right)}{\vert x\vert/(m+1)}\right)^{m+1}
     \ea
     belongs to $E_{1,m}\cap L_p(\R^m)$. Then condition (C5) is satisfied
     for $1\le p\le\iy$.

     Thus we can use Theorem \ref{T2.1} for $1\le p<\iy$ and use
     Theorem \ref{T2.2} for $p=\iy$. Finally, taking account of Proposition
     \ref{P2.4} (c) and \eqref{E2.23}, we obtain for $1\le p\le\iy$

     \ba
     \sup_{f\in \left(E_{1,m}\cap L_p(\R^m)\right)\setminus\{0\}}
     \frac{\left\vert\Delta^N (f)(0)\right\vert}
     {\left\|f\right\|_{L_p(\R^m)}}
     =\left(\frac{2}{\left\vert S^{m-1}\right\vert_{m-1}}\right)^{1/p}
     \sup_{g\in  \left(E_{1,m,e}\cap L_{p,\vert t\vert^{m-1}}(\R^1)\right)
     \setminus\{0\}}
     \frac{\left\vert \left(Be_{m/2-1}\right)^N (g)(0)\right\vert}
     {\left\|g\right\|_{L_{p,\vert t\vert^{m-1}}(\R^1)}}.
     \ea
     Thus we arrive at the following corollary.
     \begin{corollary}\label{C2.13}
     For $m\in\N,\,N\in\Z^1_+$, and $p\in[1,\iy]$,
     \bna\label{E2.24}
    && \C_0\left(E_{1,m}\cap L_p(\R^m) ,\,L_p(\R^m),\,
     \Delta^N\right)\nonumber\\
     &&=\left(\frac{2}{\left\vert S^{m-1}\right\vert_{m-1}}\right)^{1/p}
     \C_0\left(E_{1,1,e}\cap L_{p,\vert t\vert^{m-1}}(\R^1),
     \,L_{p,\vert t\vert^{m-1}}(\R^1),\,
     \left(Be_{m/2-1}\right)^N\right).
     \ena
     \end{corollary}
     \end{example}

     \begin{example}\label{Ex2.14}
     $\Om=S^{m-1},\,F(\Om)=L_p(S^{m-1}),\,m\ge 2,\,1\le p\le\iy;\,
     G_m=G_m(a)=O(m,a)$, where $a\in S^{m-1}$  is a fixed point;
     $B=\PP_{n,m}\left\vert_{S^{m-1}}\right.;\,
     D=\Delta^N\left\vert_{S^{m-1}}\right.=\boldsymbol{\de}_x^N$
     (see Proposition
     \ref{P2.4} (a)), $N\in\Z^1_+$.

     Conditions (C1), (C2), and (C3) are satisfied by Proposition
     \ref{P2.6} and (C4)  is satisfied by Proposition
     \ref{P2.7} (b).

     In addition, for $1\le p\le\iy$,
     \beq\label{E2.25}
     \PP_{n,m}\cap L_p(S^{m-1})^{O(m,a)}
     =\left\{P(x)=Q((x,a)):
     Q\in \PP_{n,1}\right\}
     \eeq
     (see \cite[Proposition 5.2]{G2008}).
     In particular, condition (C5) is satisfied by \eqref{E2.25}.

     Finally, taking account of the formula ($\vphi\in L_p(S^{m-1})$)
     \ba
     \left(\int_{S^{m-1}}\vert\vphi((x,a))\vert^pdS(x)\right)^{1/p}
     =
     \left(\left\vert S^{m-2}\right\vert_{m-2}\int_{-1}^1\vert\vphi(t)\vert^p(1-t^2)^{(m-3)/2}dt\right)^{1/p},
     \qquad m\ge 2,
     \ea
     (see, e.g., \cite[Sect. 11.4]{Erd1953}) and using \eqref{E2.25} and Proposition
     \ref{P2.4} (d), we obtain from Theorem \ref{T2.1}
     \ba
     \sup_{P\in\PP_{n,m}\setminus\{0\}}
     \frac{\left\vert\boldsymbol{\de}_x^N (P)(a)\right\vert}
     {\left\|P\right\|_{L_p(S^{m-1})}}
     =\left(\frac{1}{\left\vert S^{m-2}\right\vert_{m-2}}\right)^{1/p}
     \sup_{Q\in\PP_{n,1}
     \setminus\{0\}}
     \frac{\left\vert \left(Ge_{m/2-1}\right)^N (Q)(1)\right\vert}
     {\left\|Q\right\|_{L_{p,(1-t^2)^{(m-3)/2}}([-1,1])}}.
     \ea
     Thus we arrive at the following corollary.
     \begin{corollary}\label{C2.15}
     For $a\in S^{m-1},\,m\in\{2,\,3,\,\ldots\},\,N\in\Z^1_+$, and $p\in[1,\iy]$,
     \bna\label{E2.27}
    && \C_a\left( \PP_{n,m},\,L_p(S^{m-1}),\,
     \boldsymbol{\de}_x^N\right)\nonumber\\
     &&=
     \left(\frac{1}{\left\vert S^{m-2}\right\vert_{m-2}}\right)^{1/p}
     \C_1\left(\PP_{n,1}
     ,\,L_{p,(1-t^2)^{(m-3)/2}}([-1,1]),\,
     \left(Ge_{m/2-1}\right)^N\right).
     \ena
     \end{corollary}
     \end{example}
     Note that a similar result for $m\ge 3$ and $N=0$ was proved by Arestov and
     Deikalova \linebreak\cite[Theorem 2]{AD2013}.

     \section{Univariate and Multivariate Constants}\label{S4}
     \noindent
     \setcounter{equation}{0}
     Here, we discuss main results on limit relations between
     sharp constants in the univariate and multivariate
     Markov-Bernstein-Nikolskii type inequalities.
     The necessary notation is introduced in Sections \ref{S1}, \ref{S2},
     and \ref{S3}. In particular, the sharp constant $\C_a$ was defined by
     \eqref{E2.0} and $Be_\nu$ and $Ge_\la$ are Bessel and Gegenbauer
     operators defined by
     \eqref{E2.11} and \eqref{E2.12}, respectively.

     We first discuss sharp constants in the univariate inequalities of
     different weighted metrics.

     \begin{theorem}\label{T4.1}
     If $\nu\ge -1/2,\,N\in\Z^1_+,$ and $p\in(0,\iy]$, then
     the limit relation
     \bna\label{E4.1}
     &&\lim_{n\to\iy} n^{-2N-(2\nu+2)/p}
     \C_0\left(\PP_{2\lfloor n/2 \rfloor, 1,e}
     ,\,L_{p,\vert t\vert^{2\nu+1}}([-1,1]),\,
     \left(Be_{\nu}\right)^N\right)\nonumber\\
     &&= \C_0\left(E_{1,1,e}\cap L_{p,\vert t\vert^{2\nu+1}}(\R^1),
     \,L_{p,\vert t\vert^{2\nu+1}}(\R^1),\,
     \left(Be_{\nu}\right)^N\right).
     \ena
     is valid.
     In addition, there exists a function
     $f_0\in \left(E_{1,1,e}\cap L_{p,\vert t\vert^{2\nu+1}}(\R^1)\right)
     \setminus\{0\}$ such that
     \bna\label{E4.4a}
     &&\lim_{n\to\iy} n^{-2N-(2\nu+2)/p}
     \C_0\left(\PP_{2\lfloor n/2 \rfloor, 1,e}
     ,\,L_{p,\vert t\vert^{2\nu+1}}([-1,1]),\,
     \left(Be_{\nu}\right)^N\right)\nonumber\\
  &&=\|\left(Be_{\nu}\right)^N(f_0)\|_{L_\iy(\R^1)}/
\|f_0\|_{L_{p,\vert t\vert^{2\nu+1}}(\R^1)}.
 \ena
     \end{theorem}
     \noindent
     Note that for $\nu=-1/2$ Theorem \ref{T4.1} in more general settings
     was proved in \cite[Theorem 1.1]{G2017}.

     \begin{remark}\label{R4.2}
     The domain $D_\nu(N)$ of the operator $\left(Be_\nu\right)^N$
     is a subset of the set $C^N(\R^1)$ of all $N$ times differentiable
      functions on $\R^1$ that consists of all $h\in C^N(\R^1)$,
      satisfying the relations
      $h^{(l)}(0)\prod_{j=1}^l (2j-2\nu+1)=0,\,0\le l\le N-1$.
      In particular, $D_\nu(0)=C^0(\R^1),\,
      D_\nu(1)=\{h\in C^1(\R^1):h^\prime(0)=0\}$.
      In addition, $\PP_{2\lfloor n/2\rfloor,1,e}\subseteq D_\nu(N)$ and
      $E_{1,1,e}\subseteq D_\nu(N)$.

      If $P\in\left(\PP_{n,1}\cap D_\nu(N)\right)\setminus\{0\}$,
      then $P^*(t):=(P(t)+P(-t))/2$ belongs to $\PP_{2\lfloor n/2\rfloor,1,e}
      \setminus\{0\}$ and for $p\in[1,\iy)$,
      \ba
      \frac{\vert \left(Be_{\nu}\right)^N(P)(0)
      \vert}{\|P\|_{L_{p,\vert t\vert^{2\nu+1}}([-1,1])}}
      \le \frac{\vert \left(Be_{\nu}\right)^N (P^*)(0)\vert}
      {\|P^*\|_{L_{p,\vert t\vert^{2\nu+1}}([-1,1])}}.
      \ea
      Similarly, if
       $f\in\left(E_{1,1}\cap L_{p,\vert t\vert^{2\nu+1}}(\R^1)
        \cap D_\nu(N)\right)\setminus\{0\}$, then
        $f^*(t):=(f(t)+f(-t))/2 $ belongs to $\left(E_{1,1,e}
        \cap L_{p,\vert t\vert^{2\nu+1}}(\R^1)\right)
      \setminus\{0\}$ and for $p\in[1,\iy)$.
      \ba
      \frac{\vert \left(Be_{\nu}\right)^N (f)(0)
      \vert}{\|f\|_{L_{p,\vert t\vert^{2\nu+1}}(\R^1)}}
      \le \frac{\vert \left(Be_{\nu}\right)^N(f^*)(0)\vert}
      {\|f^*\|_{L_{p,\vert t\vert^{2\nu+1}}(\R^1)}}.
      \ea
      Therefore, it is possible to replace
      $\PP_{2\lfloor n/2\rfloor,1,e}$ by $\PP_{n,1}\cap D_\nu(N)$ and
      replace
      $E_{1,1,e}$ by $E_{1,1}\cap D_\nu(N)$ in \eqref{E4.1}
      for $p\in[1,\iy)$.
      \end{remark}

      \begin{theorem}\label{T4.3}
     If $\nu\ge -1/2,\,N\in\Z^1_+,$ and $p\in(0,\iy]$, then
     the limit relation
     \bna\label{E4.5}
     &&\lim_{n\to\iy} n^{-2N-(2\nu+2)/p}
     \C_1\left(\PP_{n,1}
     ,\,L_{p,(1-t^2)^\nu}([-1,1]),\,
     \left(Ge_{\nu+1/2}\right)^N\right)\nonumber\\
     &&= 2^{1/p}\C_0\left(E_{1,1,e}\cap L_{p,\vert t\vert^{2\nu+1}}(\R^1),
     \,L_{p,\vert t\vert^{2\nu+1}}(\R^1),\,
     \left(Be_{\nu}\right)^N\right).
     \ena
     is valid.
     In addition, there exists a function
     $f_0\in \left(E_{1,1,e}\cap L_{p,\vert t\vert^{2\nu+1}}(\R^1)\right)
     \setminus\{0\}$ such that
     \bna\label{E4.5a}
      &&\lim_{n\to\iy} n^{-2N-(2\nu+2)/p}
     \C_1\left(\PP_{n,1}
     ,\,L_{p,(1-t^2)^\nu}([-1,1]),\,
     \left(Ge_{\nu+1/2}\right)^N\right)\nonumber\\
     &&=2^{1/p}\frac{\vert \left(Be_{\nu}\right)^N (f_0)(0)
      \vert}{\|f_0\|_{L_{p,\vert t\vert^{2\nu+1}}(\R^1)}}.
     \ena
\end{theorem}
\noindent
Note that for $N=0$ and $\nu=-1/2$ relation \eqref{E4.5}
     was proved in \cite[p. 246]{LL2015a} by a different method,
     while for $\nu=-1/2$ Theorem \ref{T4.3} in more general settings
      was proved in \cite[Theorem 1.5]{GT2017}.

The proofs of Theorems \ref{T4.1} and \ref{T4.3} are presented
     in Section \ref{S6}.

     We also define
     the following sharp constant which is similar to $\C_a$:
     \ba
\C= \C\left( \B, F(\Om), D\right):=\sup_{f\in \B\setminus \{0\}}
 \frac{\left\|D(f)\right\|_{L_\iy(\Om)}}
 {\|f\|_{F(\Om)}}.
 \ea
    Certainly $\C_a\le \C$ but in some cases $\C_a=\C$.

    Next, we discuss limit relations between sharp constants
    in multivariate inequalities of different metrics.

    \begin{corollary}\label{C4.4}
    If $m\in\N,\,N\in\Z^1_+$, and $p\in[1,\iy]$, then
  \bna\label{E4.6}
&&\lim_{n\to\iy} n^{-2N-m/p}
 \C_0\left( \PP_{n,m},\,L_p(\BB^m),\,
     \Delta^N\right)\nonumber\\
   && =\C_0\left(E_{1,m}\cap L_p(\R^m) ,\,L_p(\R^m),\,
     \Delta^N\right)\nonumber\\
   && = \C\left(E_{1,m}\cap L_p(\R^m) ,\,L_p(\R^m),\,
     \Delta^N\right).
     \ena
     \end{corollary}
   \proof
   The first relation in \eqref{E4.6} follows from
   Corollaries
   \ref{C2.11} and \ref{C2.13} and Theorem \ref{T4.1}
   for $\nu=m/2-1$.
   Next, let $f\in E_{1,m}\cap L_p(\R^m),\,p\in[1,\iy]$.
   Then inequality \eqref{E2.9a} holds for $\sa=1$;
   therefore, $f$ has exponential type $1$ by the definition in
   \cite[Sect. 3.1]{N1969}. In addition,
   $\left(\Delta^N\right)(f)$ has exponential type $1$ as well
   (see \cite[Sect. 3.1]{N1969} or \cite[Lemma 2.1 (d)]{G2018}),
   and also
   $\left(\Delta^N\right)(f)\in L_p(\R^m)$ by
   Bernstein's inequality (see, e. g., \cite[Eq. 3.2.2(8)]{N1969}).

   If $p\in[1,\iy)$,
   then $\lim_{\vert x\vert\to \iy}\left(\Delta^N\right)(f)(x)=0$
   by \cite[Theorem 3.2.5]{N1969}.
   Therefore, there exists $x_0\in\R^m$ such that
   $\left\|\left(\Delta^N\right)(f)\right\|_{L_\iy(\R^m)}
   =\left\vert\left(\Delta^N\right)(f)(x_0)\right\vert$.
    Setting now $f^*(x):=f(x_0-x)$, we see that
    $\left\|\left(\Delta^N\right)(f)\right\|_{L_\iy(\R^m)}
    =\left\vert\left(\Delta^N\right)(f^*)(0)\right\vert$
    and
    $\left\|f\right\|_{L_p(\R^m)}
    =\left\|f^*\right\|_{L_p(\R^m)}$.
    Thus $\C\le \C_0$ for $p\in[1,\iy)$.

    If $p=\iy$, then for any $\vep>0$ there exists
    $x_0=x_0(\vep)\in \R^m $ such that
     $\left\|\left(\Delta^N\right)(f)\right\|_{L_\iy(\R^m)}
   <\left\vert\left(\Delta^N\right)(f)(x_0)\right\vert
   +\vep \left\|f\right\|_{L_\iy(\R^m)}$.
   Setting again
   $f^*(x):=f(a-x)$, we see that
    $\left\|\left(\Delta^N\right)(f)\right\|_{L_\iy(\R^m)}
    <\left\vert\left(\Delta^N\right)(f^*)(0)\right\vert
    +\vep \left\|f\right\|_{L_\iy(\R^m)}$
    and
    $\left\|f\right\|_{L_\iy(\R^m)}
    =\left\|f^*\right\|_{L_\iy(\R^m)}$.
    Thus $\C\le \C_0+\vep$ for $p=\iy$,
    and the second equality in \eqref{E4.6}
    is established.\hfill $\Box$

    \begin{corollary}\label{C4.5}
    If $m\in\{2,\,3,\,\ldots\},\,N\in\Z^1_+,\,p\in[1,\iy],$
    and $a\in S^{m-1}$, then
  \bna\label{E4.7}
&&\lim_{n\to\iy} n^{-2N-(m-1)/p}
 \C\left( \PP_{n,m},\,L_p(S^{m-1}),\,
     \boldsymbol{\de}_x^N\right)\nonumber\\
    && =\lim_{n\to\iy} n^{-2N-(m-1)/p}
 \C_a\left( \PP_{n,m},\,L_p(S^{m-1}),\,
     \boldsymbol{\de}_x^N\right)\nonumber\\
   && =\C_0\left(E_{1,m}\cap L_p(\R^m) ,\,L_p(\R^m),\,
     \Delta^N\right)\nonumber\\
   && = \C\left(E_{1,m}\cap L_p(\R^m) ,\,L_p(\R^m),\,
     \Delta^N\right).
     \ena
     \end{corollary}
   \proof
The third relation in \eqref{E4.7} is
proved in the proof of Corollary \ref{C4.4}
and the first one can be proved similarly.
 Finally, the second equality
 follows from Corollaries
   \ref{C2.13} and \ref{C2.15} and Theorem \ref{T4.3}
   for $\nu=(m-3)/2$.\hfill $\Box$

   Note that the following special case of \eqref{E4.7} for $N=0$,
   \ba
   \lim_{n\to\iy} n^{-(m-1)/p}
 \C\left( \PP_{n,m},\,L_p(S^{m-1}),\,I\right)
  =\C\left(E_{1,m}\cap L_p(\R^m) ,\,L_p(\R^m),\,I\right),
  \ea
   was proved in
   \cite[Theorem 1.1 (i)]{DGT2018} by a different method.
   The authors of \cite{DGT2018} state that their proof of
   Theorem 1.1 is fairly nontrivial compared with
   \cite{G2005, LL2015a, LL2015b, GT2017}. In this paper we
   show that an approach to limit relations between sharp
   constants developed in \cite{GT2017} can be applied
   to even more general relations than those in \cite{DGT2018}.
   A general approach to these problems is developed in \cite{G2019}.

   Finally, we discuss an asymptotic relation between sharp constants
    in the classical univariate Nikolskii-type inequality.

    \begin{corollary}\label{C4.6}
    If $p\in[1,\iy)$, then
    \bna\label{E4.8}
    \lim_{n\to\iy} n^{-2/p}
 \C\left( \PP_{n,1},\,L_p([-1,1]),\,I\right)
  =2^{1/p}\C_0\left(E_{1,1}
  \cap L_{p,\vert t\vert}(\R^1),\,L_{p,\vert t\vert}(\R^1),\,I\right).
    \ena
    \end{corollary}
    \proof
    We first note that
     \bna\label{E4.9}
     \C_0\left(E_{1,1}
  \cap L_{p,\vert t\vert}(\R^1),\,L_{p,\vert t\vert}(\R^1),\,I\right)
  =\C_0\left(E_{1,1,e}
  \cap L_{p,\vert t\vert}(\R^1),\,L_{p,\vert t\vert}(\R^1),\,I\right).
     \ena
     This equality follows from Theorem \ref{T2.1} since conditions
     (C1) through (C5) are obviously satisfied for
     $m=1,\,B=E_{1,1}
  \cap L_{p,\vert t\vert}(\R^1),\,G_1=\{-e^*,e^*\},\,D=I,$ and
  $F(\Omega_1)=L_{p,\vert t\vert}(\R^1)$,
  where $e^*$ is the identity transformation on $\R^1$.
  In other words, \eqref{E4.9} follows from the following simple fact:
  if $f\in E_{1,1}
  \cap L_{p,\vert t\vert}(\R^1)$, then
  $(f(\cdot)+f(-\cdot))/2 \in E_{1,1}
  \cap L_{p,\vert t\vert}(\R^1)$.

  Next, Arestov and Deikalova \cite[Theorem 1]{AD2015} proved that
  \bna\label{E4.10}
   \C\left( \PP_{n,1},\,L_p([-1,1]),\,I\right)
   = \C_1\left( \PP_{n,1},\,L_p([-1,1]),\,I\right).
  \ena
   Then   \eqref{E4.8} follows from
   equalities \eqref{E4.9} and \eqref{E4.10} and
   relation \eqref{E4.5} for $N=0$
   and $\nu=0$.
   \hfill $\Box$

   Note that a different asymptotic relation for
   $\C\left( \PP_{n,1},\,L_p([-1,1]),\,I\right),\,p\in(0,\iy],$
   was proved in \cite[Theorem 1.4]{G2017}
   (see also \cite[p. 94]{G2017}).

  \section{Properties of Entire Functions and
  Polynomials}\label{S5}
 \noindent
\setcounter{equation}{0}
Throughout the section we use the notation $D_1(f)(t):=f^\prime(t)/t$.
In this section we discuss certain properties of univariate entire functions
 of exponential type and polynomials
that are needed for the proofs of Theorems \ref {T4.1} and \ref{T4.3}.
We start with estimates of the error of polynomial approximation
  for functions from $E_{1,1,e}$ and Bernstein- and Nikolskii-type inequalities.

\begin{lemma}\label{L5.1}
  For $\tau\in(0,1)$ and any function $f\in E_{1,1,e}\cap L_\iy(\R^1),$
   there is a sequence of polynomials
  $P_k\in\PP_{2\lfloor k/2 \rfloor, 1,e},\,k\in\N$, such that
  \beq\label{E5.1}
  \|f-P_k\|_{L_\iy([-k\tau,k\tau])}
  \le 2(1-\tau)^{-1/2}\exp(-C_1 k)\|f\|_{L_\iy(\R^1)},\qquad k\in\N,
  \eeq
  where $C_1:=(2/3)(1-\tau)^{3/2}$.
  \end{lemma}
  This result was proved by Bernstein \cite{B1946}
  (see also \cite[Sect. 5.4.4]{T1963} and
  \cite[Appendix, Sect. 83]{A1965}). More precise and more general
  inequalities were obtained by the author
  in \cite{G1982}
  and \cite{G1991}.

  \begin{lemma}\label{L5.2}
  For  $\tau\in(0,1)$ and any $f\in  E_{1,1,e}\cap L_\iy(\R^1)$,
  there is a sequence of polynomials
  $P_n\in\PP_{2\lfloor n/2 \rfloor, 1,e}, n\in\N$, such that for $l\in\Z_+^1$,
    $r\in(0,\iy]$,
   $\al > -1$, and $\be>-1$,
  \bna
  && \lim_{n\to\iy}\left\|f^{(l)}-P_n^{(l)}
  \right\|_{L_{r,\vert t\vert^{\al}
  \left(1-(t/(\tau n))^2\right)^\be}([-\tau n,\tau n])}=0,\label{E5.2a}\\
  &&\lim_{n\to\iy}\left\|\left(Be_\nu\right)^l (f)-\left(Be_\nu\right)^l (P_n)
  \right\|_{L_{r,\vert t\vert^{\al}
  \left(1-(t/(\tau n))^2\right)^\be}([-\tau n,\tau n])}=0.\label{E5.2}
  \ena
  \end{lemma}
  \proof
  First of all, for $P_k\in\PP_{2\lfloor k/2 \rfloor, 1,e},\,k\in\N$,
   and $l\in\Z_+^1$
   we need the following crude Markov-type inequalities:
  \bna
  &&\left\|P_k^{(l)}\right\|_{L_\iy([-a,a])}
  \le \left(k^2/a\right)^l\|P_k\|_{L_\iy([-a,a])},\label{E5.3aa}\\
  &&\left\|(D_1)^l (P_k)\right\|_{L_\iy([-a,a])}
  \le \left(k^2/a\right)^{2l}\|P_k\|_{L_\iy([-a,a])},
  \label{E5.3bb}\\
  &&\left\|\left(Be_\nu\right)^l (P_k)\right\|_{L_\iy([-a,a])}
  \le \left((2\nu+2)k^4/a^2\right)^l\|P_k\|_{L_\iy([-a,a])}.
  \label{E5.3}
  \ena
  Inequality \eqref{E5.3aa} follows from A. A. Markov's
  inequality \cite[Sect. 4.8.7]{T1963}, while \eqref{E5.3bb} is a consequence
  of the Mean Value Theorem and \eqref{E5.3aa}. Combining \eqref{E5.3aa} with
  \eqref{E5.3bb}, we obtain \eqref{E5.3} since
  $Be_\nu(f)(t)=f^{\prime\prime}(t)+(2\nu+1)D_1(f)(t)$.

  Next, let $\{P_k\}_{k=1}^\iy$ be the sequence of
  polynomials from Lemma \ref{L5.1}. Then using \eqref{E5.3}
  and estimate \eqref{E5.1}, we obtain
  \ba
  &&\left\|\left(Be_\nu\right)^l (f)-\left(Be_\nu\right)^l (P_n)
  \right\|_{L_\iy([-\tau n,\tau n])}
  \le \sum_{k=n}^\iy \left\|\left(Be_\nu\right)^l (P_k-P_{k+1})
  \right\|_{L_\iy([-\tau n,\tau n])}\\
  &&\le (2\nu+2)^l(\tau n)^{-2l} \sum_{k=n}^\iy (k+1)^{4l}
      \left\|P_k-P_{k+1}\right\|_{L_\iy([-\tau n,\tau n])}\\
  &&\le  (2\nu+2)^l(\tau n)^{-2l} \sum_{k=n}^\iy (k+1)^{4l}
      \left(\left\|f-P_{k}\right\|_{L_\iy([-\tau k,\tau k])}
 + \left\|f-P_{k+1}\right\|_{L_\iy([-\tau(k+1),\tau(k+1)])}\right)\\
  &&\le  4(2\nu+2)^l(\tau n)^{-2l}(1-\tau)^{-1/2}
  \sum_{k=n}^\iy (k+1)^{4l}\exp(-C_1 k)\,\|f\|_{L_\iy(\R^1)}.
  \ea
  Hence for $n\in\N,\,l\in\Z_+^1,\,
    r\in(0,\iy],\,
   \al > -1$, and $\be>-1$, we have
  \ba
  &&\left\|\left(Be_\nu\right)^l (f)-\left(Be_\nu\right)^l (P_n)
  \right\|_{L_{r,\vert t\vert^{\al}
  \left(1-(t/(\tau n))^2\right)^\be}([-\tau n,\tau n])}\\
  &&\le \left( \mbox{B}((\al+1)/2,\,\be+1)(\tau n)^{\al+1}\right)^{1/r}
   \left\|\left(Be_\nu\right)^l (f)
  -\left(Be_\nu\right)^l (P_n)
  \right\|_{L_\iy([-\tau n,\tau n])}\\
  &&\le C_2(\tau,r,\al,\be,l) n^{(\al+1)/r-2l}
      \int_{n}^\iy y^{4l}\exp(-C_1 y)dy\, \|f\|_{L_\iy(\R^1)}\\
   &&\le C_3(\tau,r,\al,\be,l) n^{(\al+1)/r+2l}\exp(-C_1 n)\,\|f\|_{L_\iy(\R^1)}.
  \ea
  Thus \eqref{E5.2} is established.
  Relation \eqref{E5.2a}  can be proved similarly if
  we use Lemma \ref{L5.1} and inequality \eqref{E5.3aa}.
  \hfill $\Box$

  \begin{lemma}\label{L5.2a}
  (a)
  If $l\in\Z_+^1$, and $f\in  E_{\sa,1,e}\cap L_{\iy}(\R^1)$, then
  the following Bernstein-type inequalities hold:
  \beq\label{E5.3a}
  \left\| f^{(l)}\right\|_{L_{\iy}(\R^1)}
  \le \sa^l\left\|f\right\|_{L_{\iy}(\R^1)},\qquad
  \left\| (D_1)^l(f)\right\|_{L_{\iy}(\R^1)}
  \le  \sa^{2l} \left\|f\right\|_{L_{\iy}(\R^1)}.
  \eeq
  (b) Let $\al\ge 0$ and $p\in(0,\iy)$.
  If $f\in  E_{1,1,e}\cap  L_{p,\vert t\vert^{\al}}(\R^1)$, then
  the following Nikolskii-type inequality holds:
  \beq\label{E5.3b}
  \left\|f\right\|_{L_{\iy}(\R^1)}
  \le C(p,\al)
  \left\|f\right\|_{L_{p,\vert t\vert^{\al}}(\R^1)}.
  \eeq
  \end{lemma}
  \proof
 (a)  The first inequality in \eqref{E5.3a}
  is a classical Bernstein inequality \cite[Sect. 4.8.2]{T1963} and the second one
  immediately follows from the Mean Value Theorem
  and the first one.\\
  (b) Inequality \eqref{E5.3b} for $p\in[1,\iy)$ and $\al\ge 0$
   was proved by Platonov
  \cite[Theorem 3.5]{P2007}, while for $p\in(0,\iy)$ and $\al= 0$,
   \eqref{E5.3b} follows from Nikolskii's inequality \cite[Theorem 2.3.5]{N1969}.
   If $p\in(0,1)$ and $f\in  E_{1,1,e}\cap  L_{p,\vert t\vert^{\al}}(\R^1)$,
   then $f\in L_{p}(\R^1)$; hence $f\in L_{\iy}(\R^1)$ by
    Nikolskii's inequality.
   Then $f\in  L_{1,\vert t\vert^{\al}}(\R^1)$ since
   \beq\label{E5.3c}
   \left\|f\right\|_{L_{1,\vert t\vert^{\al}}(\R^1)}
   \le \left\|f\right\|_{L_{\iy}(\R^1)}^{1-p}
   \left\|f\right\|_{L_{p,\vert t\vert^{\al}}(\R^1)}^p.
   \eeq
   Using estimate \eqref{E5.3c} and Platonov's inequality \eqref{E5.3b} for $p=1$,
    we obtain
   \beq\label{E5.3d}
   \left\|f\right\|_{L_{\iy}(\R^1)}
   \le C(1,\al)  \left\|f\right\|_{L_{1,\vert t\vert^{\al}}(\R^1)}
   \le C(1,\al) \left\|f\right\|_{L_{\iy}(\R^1)}^{1-p}
   \left\|f\right\|_{L_{p,\vert t\vert^{\al}}(\R^1)}^p.
   \eeq
   Therefore, \eqref{E5.3b} for $p\in(0,1)$ follows from \eqref{E5.3d}
   with $C(p,\al)\le (C(1,\al))^{1/p}$.
\hfill $\Box$

  In addition to a compactness theorem for entire functions of exponential type from
  Proposition \ref{P2.3}, we need a different type of a compactness theorem.
  \begin{lemma}\label{L5.3}
     Let $\E_1$ be the set  of all univariate entire functions
     $f(z)=\sum_{k=0}^\iy c_kz^{2k}$,
     satisfying the following condition: for any $\de>0$ there exists
     a constant $C(\de)$,
     independent of $f$ and $k$, such that
     \beq\label{E5.3e}
     \vert c_k\vert \le \frac{C(\de)(1+\de)^{2k}}{(2k)!},\qquad k\in\Z_+^1.
     \eeq
     Then for any sequence $\{f_n\}_{n=1}^\iy\subseteq\E_1$ there exist
     a subsequence $\{f_{n_s}\}_{s=1}^\iy$ and a function $f_0\in E_{1,1,e}$
     such that for every $l\in\Z_+^1$,
      \beq\label{E5.3f}
     \lim_{s\to\iy}f_{n_s}^{(l)}=f_0^{(l)},\qquad
     \lim_{s\to\iy}(Be_\nu)^l(f_{n_s})=(Be_\nu)^l(f_{0}),
     \eeq
     uniformly on each compact subset of $\CC$.
     \end{lemma}
     \proof
     The existence of a subsequence
     $\{f_{n_s}(z)=\sum_{k=0}^\iy c_{k,n_s}z^{2k}\}_{s=1}^\iy$
     and a function \linebreak
     $f_0(z)=\sum_{k=0}^\iy c_{k,0}z^{2k}$ such that
     for every $l\in\Z_+^1,\,
     \lim_{s\to\iy}f_{n_s}^{(l)}=f_0^{(l)}$
     uniformly on the disk
     $d_R:=\{z\in\C:\vert z\vert\le R\},\,R>0$, was proved in
  \cite[Lemma 2.6 (a)]{G2017}.
  In particular,
  \beq\label{E5.3g}
  \lim_{s\to\iy} c_{k,n_s}=c_{k,0},\qquad k\in\Z_+^1.
  \eeq
  Next, it is easy to prove by induction in $l$ that if
  $f(z)=\sum_{p=0}^\iy c_p\,z^{2p}$, then
  \bna
 (Be_\nu)^l(f)(z)
  =2^{2l}\sum_{p=0}^\iy \prod_{d=1}^l(p+d)(p+d+\nu)
  \,c_{p+l}\,z^{2p}.\label{E5.3i}
  \ena
  Then we obtain from \eqref{E5.3i} for $M\in\N$
  \ba
  &&\max_{z\in d_R}\left\vert
  (Be_\nu)^l(f_{0})(z)-(Be_\nu)^l(f_{n_s})(z)\right\vert\\
  &&\le 2^{2l}\max_{z\in d_R}
  \sum_{p=0}^{M-1} \prod_{d=1}^l(p+d)(p+d+\nu)\,
  \left\vert c_{p+l,0}-c_{p+l,n_s}\right\vert \vert z\vert^{2p}\\
  &&+ 2^{2l}\max_{z\in d_R}
  \sum_{p=M}^{\iy} \prod_{d=1}^l(p+d)(p+d+\nu)\,
  \left(\left\vert
  c_{p+l,0}\right\vert+\left\vert c_{p+l,n_s}\right\vert\right) \vert z\vert^{2p}
  =S_1+S_2,
  \ea
   where by  \eqref{E5.3e} for $\de=1$,
   \ba
   S_2\le 2^{4l+1}C(1)\sum_{p=M}^{\iy}\frac{\prod_{d=1}^l(p+d)(p+d+\nu)
   \,(2R)^{2p}}{(2p+2l)!}.
   \ea
   Further, given $\vep>0$ and $R>0$, we can choose $M=M(\vep,R)$
   such that $S_2<\vep/2$. Finally, by \eqref{E5.3g}, we can
   choose $s_0=s_0(\vep,R)\in\N$ such that $S_1<\vep/2$ for all $s\ge s_0$.
   Thus the second relation in \eqref{E5.3f} holds
   uniformly on $d_R$ as well.
   \hfill $\Box$

     Certain inequalities of different weighted metrics for univariate polynomials
     are discussed in the following lemma.
  \begin{lemma}\label{L5.4}
  For $P\in\PP_{n,1},\,k\in\Z^1_+,\,\vep\in(0,1/2),\,p\in(0,\iy),
  \,\al\ge 0,$ and $\be > -1,$ the following inequalities hold:
  \bna
 && \|P\|_{L_{\iy}([-1+\vep,1-\vep])}
  \le C_4(\al,\be,p,\vep)\,n^{(\al+1)/p}\|P\|_{L_{p,\vert t\vert^\al(1-t^2)^\be}([-1,1])},
  \label{E5.5}\\
  &&\left\vert P^{(k)}(0)\right\vert
  \le C_4(\al,\be,p,\vep)(1-\vep)^{-k}n^{k+(\al+1)/p}
  \|P\|_{L_{p,\vert t\vert^\al(1-t^2)^\be}([-1,1])}.\label{E5.6}
  \ena
  \end{lemma}
  \proof
  Inequality \eqref{E5.5} for $\al=\be=0$ and $p\in(0,\iy)$  was proved in
  \cite[Eq. (2.10)]{G2017} by using an extension of Bari's inequality
  \cite[Theorem 6]{Ba1954}
  to $p\in(0,\iy)$ (see \cite[Lemma 2.4]{G2017}). Then \eqref{E5.5}
  for $\al=\be=0$ implies the estimate
  \beq\label{E5.7}
  \|P\|_{L_{\iy}([-1+\vep,1-\vep])}
  \le C_5(p,\vep)\,n^{1/p}\|P\|_{L_{p}([-1+\vep/2,1-\vep/2])},
  \qquad p\in(0,\iy).
  \eeq
  Therefore,  \eqref{E5.5} follows from \eqref{E5.7}
  and the inequalities
  \beq\label{E5.8}
  \|P\|_{L_{p}([-1+\vep/2,1-\vep/2])}
  \le C_6\,n^{\al/p}\|P\|_{L_{p,\vert t\vert^\al}([-1+\vep/2,1-\vep/2])}
  \le C_7\,n^{\al/p}\|P\|_{L_{p,\vert t\vert^\al(1-t^2)^\be}([-1,1])},
  \eeq
  where $C_6$ and $C_7$ in \eqref{E5.8}
   are independent of $P$ and $n$.
To prove the first inequality in \eqref{E5.8}, we observe that
  \beq\label{E5.9}
  \|P\|_{L_{p,\vert t\vert^\al}([-1+\vep/2,1-\vep/2])}
  \ge \|P\|_{L_{p,\vert t\vert^\al}(\{t:1-\vep/2\ge \vert t\vert\ge C/n\})}
  \ge (C/n)^{\al/p}\|P\|_{L_{p}(\{t:1-\vep/2\ge \vert t\vert\ge C/n\})},
  \eeq
  where $C\in(0,1/3)$ is a fixed number. Next, we note that $0<C/n<1/3$, so
  by \eqref{E5.7},
  \beq\label{E5.10}
  \|P\|_{L_{p}([-C/n,C/n])}
  \le \left(2C/n\right)^{1/p}\|P\|_{L_{\iy}([-C/n,C/n])}
  \le \left(2C\right)^{1/p}C_5
  \|P\|_{L_{p}([-1+\vep/2,1-\vep/2])}.
  \eeq
  Choosing now $C:=\min\{1/3,C_5^{-p}2^{-p-1}\}$,
  we obtain from \eqref{E5.10}
  \bna\label{E5.11}
 && \|P\|_{L_{p}(\{t:1-\vep/2\ge \vert t\vert\ge C/n\})}^p
  =\|P\|_{L_{p}([-1+\vep/2,1-\vep/2])}^p
  -\|P\|_{L_{p}([-C/n,C/n])}^p\nonumber\\
  &&\ge \left(1-(2C)C_5^p\right)\|P\|_{L_{p}([-1+\vep/2,1-\vep/2])}^p
  \ge (1/2)^p\|P\|_{L_{p}([-1+\vep/2,1-\vep/2])}^p.
  \ena
  Finally, combining \eqref{E5.9} and \eqref{E5.11}, we arrive at
  the first inequality in \eqref{E5.8}. Thus \eqref{E5.5} is established.
  To prove \eqref{E5.6}, we use the estimate
  \ba
  \left\vert P^{(k)}(0)\right\vert
  \le \left(\frac{n}{1-\vep}\right)^k \|P\|_{L_{\iy}([-1+\vep,1-\vep])}
  \ea
  (see \cite[Eq. 4.8(49)]{T1963}) and inequality \eqref{E5.5}.
  \hfill $\Box$

  In the next lemma, in particular, we discuss a relation between
  the Bessel and Gegenbauer operators.

  \begin{lemma}\label{L5.5}
  (a) If $b\ne 0$, then
  \beq\label{E5.12}
  \PP_{2n,1,e}
  =S:=\{P_{2n}(t)=R_n\left(1-2b^{-2}t^2\right):R_n\in\PP_{n,1}\}.
  \eeq
  (b) If  $R_n\in\PP_{n,1},\,b\ne 0,\,\nu\ge -1/2$, and $N\in\Z^1_+$, then
  \beq\label{E5.13}
  \left(Ge_{\nu+1/2}\right)^N(R_n)\left(1-2b^{-2}t^2\right)
  =\left(D_{\nu,b}\right)^N\left(P_{2n}\right)(t),
  \eeq
  where
  \beq\label{E5.14}
  D_{\nu,b}(g)(t)
  := \frac{b^2}{4}\left(\left(1-\frac{t^2}{b^2}\right)g^{\prime\prime}(t)
  +\left((2\nu+1)-\frac{4\nu+3}{b^2}t^2\right)\frac{g^\prime(t)}{t}\right),
  \eeq
  and
  \beq\label{E5.14a}
  P_{2n}(t)=P_{2n,b}(t):=R_n\left(1-2b^{-2}t^2\right).
  \eeq
  (c)
  Let $ P_{2n}=P_{2n,b(n)}$ be defined by \eqref{E5.14a},
   where $R_n\in\PP_{n,1}$
   and $b=b(n)$
  satisfies the condition $\lim_{n\to\iy}b(n)=\iy$. Next,  let
   there exist a sequence of natural numbers $\{n_s\}_{s=1}^\iy$ and
  an even entire function $f$ such that for every
  $l\in\Z_+^1$,
  \beq\label{E5.15a}
  \lim_{s\to\iy}
  P_{2n_s,b(n_s)}^{(l)}=f^{(l)},\qquad
  \lim_{s\to\iy}\left(Be_\nu\right)^l
  \left(P_{2n_s,b(n_s)}\right)=\left(Be_\nu\right)^l (f),
  \eeq
     uniformly on each compact subset of $\R^1$. Then the
     following relation holds
     for each $t\in\R^1$:
    \beq\label{E5.15}
    \lim_{s\to\iy}(b(n_s)/2)^{-2N}
    \left(Ge_{\nu+1/2}\right)^N(R_{n_s})\left(1-2(b(n_s))^{-2}t^2\right)
    =\left(Be_\nu\right)^N(f)(t).
    \eeq
  \end{lemma}
  \proof
  (a) It suffices to prove that $\PP_{2n,1,e}\subseteq S$. Let
  $P_{2n}\in \PP_{2n,1,e}$ and let
  $P_{2n}(t)=V_n(t^2)$, where $V_n\in\PP_{n,1}$.
  Then by Taylor's formula, $P_{2n}(t)=R_n\left(1-2b^{-2}t^2\right)$, where
  \ba
  R_n(y):=\sum_{k=0}^n\frac{(-1)^kV_n^{(k)}(b^2/2)b^{2k}}{k!\,2^{k}}y^k.
  \ea
  Therefore, $P_{2n}\in S$ and \eqref{E5.12} is established.\vspace{.12in}\\
  (b) Setting $D^*(f)(y):=A^*(y)f^{\prime\prime}(y)+B^*(y)f^\prime(y)$,
  we see that for $y=\vphi(t)$,
  \beq\label{E5.16}
  D^*(f)(\vphi(t))
  =\frac{A^*(\vphi(t))}{\vphi^{\prime 2}(t)} H^{\prime\prime}(t)
  +\left(\frac{B^*(\vphi(t))}{\vphi^\prime(t)}
  -\frac{A^*(\vphi(t))\vphi^{\prime\prime}(t)}
  {\vphi^{\prime 3}(t)}\right)H^{\prime}(t),
  \eeq
  where $H(t):=f(\vphi(t))$. Then choosing
  \ba
  A^*(y)=1-y^2,\quad B^*(y)=-(2\nu+2)y,\quad f(y)=R_n(y),
   \quad\vphi(t)=1-2b^{-2}t^2, \quad H(t)=P_{2n,b}(t),
  \ea
  we obtain \eqref{E5.13} for $N=1$ from \eqref{E5.16}
  by a straightforward calculation.
  Next, it follows from \eqref{E5.13} for $N=1$ that for $k\in\N$
  \ba
   \left(Ge_{\nu+1/2}\right)^k(R_n)\left(1-2b^{-2}t^2\right)
  =D_{\nu,b}\left(\left(Ge_{\nu+1/2}\right)^{k-1}(R_n)\left
  (1-2b^{-2}t^2\right)\right).
  \ea
  Hence
  identity \eqref{E5.13}
  can be proved by induction in $k$.\vspace{.12in}\\
  (c)
  By \eqref{E5.13} and \eqref{E5.14},
  \beq\label{E5.17}
  (b(n_s)/2)^{-2N}\left(D_{\nu,b(n_s)}\right)^N
  \left(P_{2n_s,b(n_s)}\right)(t)
  =\left(Be_\nu\left( P_{2n_s,b(n_s)}\right)
  -(b(n_s))^{-2}U\left( P_{2n_s,b(n_s)}\right)\right)^N(t),
  \eeq
  where for any $A>0,\,U(g)(t):=t^2g^{\prime\prime}(t)+tg^\prime(t)$
  is a continuous differential operator in the $L_\iy([-A,A])$-metric
  on the set of all even entire functions.
  Then using \eqref{E5.15a} and \eqref{E5.17}, we see that
  \beq\label{E5.18}
  \lim_{s\to\iy}(b(n_s)/2)^{-2N}\left(D_{\nu,b(n_s)}\right)^N
  \left(P_{2n_s,b(n_s)}\right)(t)
  =\left(Be_{\nu}\right)^N(f)(t)
  \eeq
  uniformly on each interval $[-A,A],\,A>0$.
  Thus \eqref{E5.15} follows from \eqref{E5.13} and \eqref{E5.18}.
  \hfill $\Box$

 \section{Proofs of Theorems  \ref{T4.1} and \ref{T4.3}}\label{S6}
 \noindent
\setcounter{equation}{0}
Throughout the section we use the notation $\tilde{p}=\min\{1,p\}$
for $p>0$ introduced in Section \ref{S3}
and also use the operator $D_1(f)(t)=f^\prime(t)/t$
 introduced in Section \ref{S5}.\vspace{.12in}\\
\emph{Proof of Theorem \ref{T4.1}.}
We first prove the inequality
\bna\label{E6.1}
&&\C_0\left(E_{1,1,e}\cap L_{p,\vert t\vert^{2\nu+1}}(\R^1),
     \,L_{p,\vert t\vert^{2\nu+1}}(\R^1),\,
     \left(Be_{\nu}\right)^N\right)\nonumber\\
&&\le \liminf_{n\to\iy} n^{-2N-(2\nu+2)/p}
     \C_0\left(\PP_{2\lfloor n/2 \rfloor, 1,e}
     ,\,L_{p,\vert t\vert^{2\nu+1}}([-1,1]),\,
     \left(Be_{\nu}\right)^N\right).
     \ena
Let $f$ be any function from
$  E_{1,1,e}\cap L_{p,\vert t\vert^{2\nu+1}}(\R^1),\,
p\in(0,\iy],\,\nu\ge -1/2$,
and let $\tau\in(0,1)$
 be a fixed number.
Then using  even polynomials
$P_n\in\PP_{2\lfloor n/2 \rfloor, 1,e}$
 from Lemma \ref{L5.2} for $r=\iy$, we obtain by \eqref{E5.2}
 and by definition \eqref{E2.0} of
 $\C_0$,
\bna\label{E6.3}
&&\left\vert \left(Be_{\nu}\right)^N(f)(0)\right\vert
 =\lim_{n\to\iy}\left\vert
 \left(Be_{\nu}\right)^N(P_n)(0)\right\vert\nonumber\\
 &&\le \liminf_{n\to\iy}\left(
 \C_0\left(\PP_{2\lfloor n/2 \rfloor, 1,e}
     ,\,L_{p,\vert t\vert^{2\nu+1}}([-n\tau,n\tau]),\,
     \left(Be_{\nu}\right)^N\right)
     \|P_n\|_{L_{p,\vert t\vert^{2\nu+1}}([-n\tau,n\tau])}\right).
 \ena
Next, note that $f\in L_\iy(\R^1)$, by Nikolskii-type
inequality \eqref{E5.3b}.
 Further, applying "triangle" inequality \eqref{E2.7a}
 and using again relation \eqref{E5.2} of
  Lemma \ref{L5.2} for $\al=2\nu+1,\, \be=0,\,l=0$,
  and $r=p,\,p\in(0,\iy]$,
 we have
 \bna\label{E6.4}
 &&\limsup_{n\to\iy}
 \left\| P_n\right\|_{L_{p,\vert t\vert^{2\nu+1}}
 ([-n\tau,n\tau])}^{\tilde{p}}\nonumber\\
 &&\le \lim_{n\to\iy}\left(\|f-P_n\|_{L_{p,\vert t\vert^{2\nu+1}}
 ([-n\tau,n\tau])}^{\tilde{p}}
 +\|f\|_{L_{p,\vert t\vert^{2\nu+1}}([-n\tau,n\tau])}^{\tilde{p}}\right)
 =\|f\|_{L_{p,\vert t\vert^{2\nu+1}}(\R^1)}^{\tilde{p}}.
 \ena
 Combining \eqref{E6.3} with \eqref{E6.4}, we obtain
 \bna\label{E6.4a}
 &&\left\vert \left(Be_{\nu}\right)^N(f)(0)\right\vert/
 \|f\|_{L_{p,\vert t\vert^{2\nu+1}}(\R^1)}\nonumber\\
 &&\le \liminf_{n\to\iy}
 \C_0\left(\PP_{2\lfloor n/2 \rfloor, 1,e}
     ,\,L_{p,\vert t\vert^{2\nu+1}}([-n\tau,n\tau]),\,
     \left(Be_{\nu}\right)^N\right)\nonumber\\
     &&=\tau^{-2N-(2\nu+2)/p} \liminf_{n\to\iy}n^{-2N-(2\nu+2)/p}
 \C_0\left(\PP_{2\lfloor n/2 \rfloor, 1,e}
     ,\,L_{p,\vert t\vert^{2\nu+1}}([-1,1]),\,
     \left(Be_{\nu}\right)^N\right).
     \ena
 Letting $\tau\to 1-$ in \eqref{E6.4a},
 we arrive at \eqref{E6.1} for $\nu\ge -1/2,\,N\in\Z_+^1$, and $p\in(0,\iy]$.

  Further, we  prove the
 inequality
 \bna\label{E6.6}
&&\limsup_{n\to\iy} n^{-2N-(2\nu+2)/p}
     \C_0\left(\PP_{2\lfloor n/2 \rfloor, 1,e}
     ,\,L_{p,\vert t\vert^{2\nu+1}}([-1,1]),\,
     \left(Be_{\nu}\right)^N\right)\nonumber\\
     &&\le \C_0\left(E_{1,1,e}\cap L_{p,\vert t\vert^{2\nu+1}}(\R^1),
     \,L_{p,\vert t\vert^{2\nu+1}}(\R^1),\,
     \left(Be_{\nu}\right)^N\right),
     \ena
  by constructing a nontrivial function $f_0\in E_{1,1,e}
  \cap L_{p,\vert t\vert^{2\nu+1}}(\R^1)$
   such that
 \bna \label{E6.7}
 &&\limsup_{n\to\iy} n^{-2N-(2\nu+2)/p}
     \C_0\left(\PP_{2\lfloor n/2 \rfloor, 1,e}
     ,\,L_{p,\vert t\vert^{2\nu+1}}([-1,1]),\,
     \left(Be_{\nu}\right)^N\right)\nonumber\\
  &&\le\left\vert \left(Be_{\nu}\right)^N(f_0)(0)\right\vert/
\|f_0\|_{L_{p,\vert t\vert^{2\nu+1}}(\R^1)}\nonumber\\
  &&\le \C_0\left(E_{1,1,e}\cap L_{p,\vert t\vert^{2\nu+1}}(\R^1),
     \,L_{p,\vert t\vert^{2\nu+1}}(\R^1),\,
     \left(Be_{\nu}\right)^N\right).
 \ena
 Then inequalities \eqref{E6.1} and \eqref{E6.6} imply  \eqref{E4.1}.
 In addition, $f_0$ is an extremal function in \eqref{E4.1},
that is, relation \eqref{E4.4a} is valid.

It remains to construct a nontrivial function $f_0$,
satisfying \eqref{E6.7}.
We first note that
\beq \label{E6.8}
\inf_{n\ge 2N+2} n^{-2N-(2\nu+2)/p}
\C_0\left(\PP_{2\lfloor n/2 \rfloor, 1,e}
     ,\,L_{p,\vert t\vert^{2\nu+1}}([-1,1]),\,
     \left(Be_{\nu}\right)^N\right)
     \ge C_{8}(p,N,\nu).
\eeq
This inequality follows immediately from  \eqref{E6.1}.
Let $P_n\in\PP_{2\lfloor n/2 \rfloor, 1,e}$
 be an even polynomial, satisfying the equality
\beq \label{E6.9}
\C_0\left(\PP_{2\lfloor n/2 \rfloor, 1,e}
     ,\,L_{p,\vert t\vert^{2\nu+1}}([-1,1]),\,
     \left(Be_{\nu}\right)^N\right)
=\frac{\left\vert \left(Be_{\nu}\right)^N(P_n)(0)\right\vert}
{\|P_n\|_{L_{p,\vert t\vert^{2\nu+1}}([-1,1])}},\quad n\in\N.
\eeq
The existence of an extremal polynomial $P_n$ in \eqref{E6.9}
can be proved by the standard compactness argument (cf. \cite{GT2017}).
Next, setting $Q_n(x):=P_n(x/n)$, we have from \eqref{E6.9} that
\bna \label{E6.10}
&&n^{-2N-(2\nu+2)/p}
\C_0\left(\PP_{2\lfloor n/2 \rfloor, 1,e}
     ,\,L_{p,\vert t\vert^{2\nu+1}}([-1,1]),\,
     \left(Be_{\nu}\right)^N\right)\nonumber\\
&&=\left\vert \left(Be_{\nu}\right)^N(Q_n)(0)\right\vert
/\|Q_n\|_{L_{p,\vert t\vert^{2\nu+1}}([-n,n])}, \qquad n\in\N.
\ena
We can assume that
\beq \label{E6.11}
  \left(Be_{\nu}\right)^N(Q_n)(0)=1,\qquad n\in\N.
\eeq
Then it follows from \eqref{E6.10}, \eqref{E6.11}, and \eqref{E6.8}
that for $n\ge 2N+2$,
\bna \label{E6.12}
\|Q_n\|_{L_{p,\vert t\vert^{2\nu+1}}([-n,n])}
&=&\left(n^{-2N-(2\nu+2)/p}
\C_0\left(\PP_{2\lfloor n/2 \rfloor, 1,e}
     ,\,L_{p,\vert t\vert^{2\nu+1}}([-1,1]),\,
     \left(Be_{\nu}\right)^N\right)\right)^{-1}\nonumber\\
&\le& 1/C_{8}(p,N,\nu),\qquad n\in\N.
\ena
Further, $Q_n\in\PP_{n,1}$ and it follows from
inequality \eqref{E5.6} of Lemma \ref{L5.4}
for $\al=2\nu+1$ and $\be=0$ and
from \eqref{E6.12}  that
for any $\vep\in(0,1/2)$ and any $k\in\Z_+^1$,
\bna \label{E6.13}
\left\vert Q_n^{(k)}(0)\right\vert
&\le& C_4(2\nu+1,0,p,\vep)
(1-\vep)^{-k }\|Q_n\|_{L_{p,\vert t\vert^{2\nu+1}}([-n,n])}\nonumber\\
&\le& \left(C_4(2\nu+1,0,p,\vep)/C_{8}(p,N,\nu)\right)(1-\vep)^{-k}.
\ena
Let $\{n_r\}_{r=1}^\iy$ be a subsequence of natural numbers such that
\bna \label{E6.14}
&&\limsup_{n\to\iy}\,
n^{-2N-(2\nu+2)/p}
\C_0\left(\PP_{2\lfloor n/2 \rfloor, 1,e}
     ,\,L_{p,\vert t\vert^{2\nu+1}}([-1,1]),\,
     \left(Be_{\nu}\right)^N\right)\nonumber\\
     &&=\lim_{r\to\iy}
     (n_r)^{-2N-(2\nu+2)/p}
\C_0\left(\PP_{2\lfloor n_r/2 \rfloor, 1,e}
     ,\,L_{p,\vert t\vert^{2\nu+1}}([-1,1]),\,
     \left(Be_{\nu}\right)^N\right).
\ena
Inequality \eqref{E6.13} shows that the polynomial sequence
$\{Q_{n_r}\}_{r=1}^\iy\subseteq \E_1$ satisfies the conditions
of Lemma \ref{L5.3}. Therefore, there exist
 a function $f_0\in E_{1,1,e}$ and a subsequence
 $\{Q_{n_{r_s}}\}_{s=1}^\iy$ such that
 \beq \label{E6.15}
\lim_{s\to\iy}(B_\nu)^l\left(Q_{n_{r_s}}\right)(t)=
 (B_\nu)^l(f_0)(t),
 \qquad 0\le l\le 2N,
 \eeq
 uniformly on any interval $[-A,A],\,A>0$. Moreover, by
  \eqref{E6.15} for $l=2N$ and \eqref{E6.11},
  \beq \label{E6.16}
  \left(Be_{\nu}\right)^N(f_0)(0)=1.
\eeq
In addition,   applying "triangle" inequality \eqref{E2.7a}
and using \eqref{E6.15} for $l=0$, \eqref{E6.10}, and \eqref{E6.11},
we obtain for any interval $[-A,A],\,A>0$,
\bna \label{E6.17}
&&\|f_0\|_{L_{p,\vert t\vert^{2\nu+1}}([-A,A])}^{\tilde{p}}\nonumber\\
&&\le \lim_{s\to\iy}
\left(\|f_0-Q_{n_{r_s}}\|_{L_{p,\vert t\vert^{2\nu+1}}([-A,A])}^{\tilde{p}}
+\|Q_{n_{r_s}}\|_{L_{p,\vert t\vert^{2\nu+1}}([-A,A])}^{\tilde{p}}\right)\nonumber\\
&&\le \lim_{s\to\iy}\|Q_{n_{r_s}}
\|_{L_{p,\vert t\vert^{2\nu+1}}([-n_{r_s},n_{r_s}])}^{\tilde{p}}
\nonumber\\
&&=
\left(\lim_{s\to\iy}
(n_{r_s})^{-2N-(2\nu+1)/p}
\C_0\left(\PP_{2\lfloor n_{r_s}/2 \rfloor, 1,e}
     ,\,L_{p,\vert t\vert^{2\nu+1}}([-1,1]),\,
     \left(Be_{\nu}\right)^N\right)\right)^{-\tilde{p}}.
\ena
Next using \eqref{E6.17} and \eqref{E6.8}, we see that
\beq \label{E6.18}
\|f_0\|_{L_{p,\vert t\vert^{2\nu+1}}(\R^1)}\le 1/C_{8}(p,N,\nu).
\eeq
Therefore, $f_0$ is a nontrivial function from
$E_{1,1.e}\cap L_{p,\vert t\vert^{2\nu+1}}(\R^1)$,
by \eqref{E6.16} and \eqref{E6.18}. Thus for any interval
$[-A,A],\,A>0$, we obtain from \eqref{E6.14}, \eqref{E6.15},
 and \eqref{E6.16}
\bna \label{E6.19}
&&\limsup_{n\to\iy} n^{-2N-(2\nu+2)/p}
     \C_0\left(\PP_{2\lfloor n/2 \rfloor, 1,e}
     ,\,L_{p,\vert t\vert^{2\nu+1}}([-1,1]),\,
     \left(Be_{\nu}\right)^N\right)\nonumber\\
&&=\lim_{s\to\iy}\left(\|Q_{n_{r_s}}\|
_{L_{p,\vert t\vert^{2\nu+1}}([-n_{r_s},n_{r_s}])}\right)^{-1}\nonumber\\
&&\le \lim_{s\to\iy}\left(\|Q_{n_{r_s}}\|
_{L_{p,\vert t\vert^{2\nu+1}}([-A,A])}\right)^{-1}\nonumber\\
&&=\left\vert \left(Be_{\nu}\right)^N (f_0)(0)\right\vert/
\|f_0\|_{L_{p,\vert t\vert^{2\nu+1}}([-A,A])}.
\ena
Finally, letting $A\to \iy$ in \eqref{E6.19}, we arrive at \eqref{E6.7}
for $\nu\ge -1/2,\,N\in\Z_+^1$, and $p\in(0,\iy]$.\hfill$\Box$
\vspace{.12in}\\
\emph{Proof of Theorem \ref{T4.3}.}
The proof is similar to the proof of Theorem \ref{T4.1} but it needs more
 technical details.
 We first prove the inequality
\bna\label{E6.20}
&&2^{1/p}\C_0\left(E_{1,1,e}\cap L_{p,\vert t\vert^{2\nu+1}}(\R^1),
     \,L_{p,\vert t\vert^{2\nu+1}}(\R^1),\,
     \left(Be_{\nu}\right)^N\right)\nonumber\\
&&\le \liminf_{n\to\iy} n^{-2N-(2\nu+2)/p}
     \C_1\left(\PP_{n,1}
     ,\,L_{p,(1-u^2)^\nu}([-1,1]),\,
     \left(Ge_{\nu+1/2}\right)^N\right).
     \ena
Let $f$ be any function from
$  E_{1,1,e}\cap L_{p,\vert t\vert^{2\nu+1}}(\R^1),\,p\in(0,\iy]$,
and let $\tau\in(0,1)$
 be a fixed number. It follows from Nikolskii-type inequality
 \eqref{E5.3b} that $f\in L_\iy(\R^1)$. Given $\vep\in(0,1/2)$, we define
 \beq\label{E6.20aa}
 F(t)=F_{\nu,\vep}(t):=
 \left\{\begin{array}{ll} f(t), &\nu\ge 0,\\
 f_\vep(t)g_\vep(t):=
 f((1-\vep)t)\left(\frac{\sin\left(\vep t/d\right)}{\vep t/d}\right)^d,
 &-1/2\le\nu<0,
 \end{array}\right.
 \eeq
 where $d:=\lfloor(2\nu+2)/p\rfloor +1$. It is easy to see that
 $F\in E_{1,1,e}\cap L_\iy(\R^1)$ and
 \beq\label{E6.20a}
 \|F\|_{L_{p,\vert t\vert^{2\nu+1}}(\R^1)}
 \le (1-\vep)^{-(2\nu+2)/p}\|f\|_{L_{p,\vert t\vert^{2\nu+1}}(\R^1)}.
 \eeq
 In addition, we prove the equality
 \beq\label{E6.20b}
 \left(Be_{\nu}\right)^N(F)(0)
 =(\Lambda_\nu(\vep))^{2N}\left( Be_{\nu}\right)^N(f)(0)
 +\eta_{\nu,N}(\vep),
 \eeq
 where $\eta_{\nu,N}(\vep)=0$ for $N=0$ or $\nu\ge 0$ and
 $\vert \eta_{\nu,N}(\vep)\vert \le C_9\vep \|f\|_{L_{\iy}(\R^1)}$
 for $N\ge 1$ and $\nu\in[-1/2,0)$. Here,
 $C_9$ is independent of $\vep$ and $f$, and
 \ba
 \Lambda_\nu(\vep):=\left\{\begin{array}{ll}1,&\nu\ge 0,\\
 1-\vep, &-1/2\le\nu< 0.\end{array}\right.
 \ea
 Equality \eqref{E6.20b} holds trivially for $N=0$ or $\nu\ge 0$.
 To prove \eqref{E6.20b} for $N\ge 1$ and
  $\nu\in[-1/2,0)$, we first need
 the Leibniz-type
 rule for the Bessel operator. Note that
 the Leibniz rule holds for the operator $D_1$, that is,
 \ba
 \left(D_1\right)^k(\vphi\cdot\psi)(t)
 =\sum_{l=0}^k\binom{k}{l}\left(D_1\right)^l(\vphi)(t)
 \left(D_1\right)^{k-l}(\psi)(t).
 \ea
 Taking also into account the Leibniz rule for the derivative $D$,
 we arrive at the following formula:
 \bna\label{E6.20b1}
 &&\left(Be_\nu\right)^N(\vphi\cdot\psi)(t)
 =\left((D)^2+(2\nu+1)D_1\right)^N(\vphi\cdot\psi)(t)\nonumber\\
 &&=\sum_{j_1+\ldots+j_N=N}
 (2\nu+1)^{j_2+j_4+\ldots}
 \sum_{m_1=0}^{2j_1}\sum_{m_2=0}^{j_2}\sum_{m_3=0}^{2j_3}
 \cdot\cdot\cdot
 \binom{2j_1}{m_1}\binom{j_2}{m_2}\binom{2j_3}{m_3}
 \cdot\cdot\cdot\nonumber\\
&& \left(D\right)^{m_1}\left(D_1\right)^{m_2}\left(D\right)^{m_3}
 \cdot\cdot\cdot
 (\psi)(t)
 \cdot
 \left(D\right)^{2j_1-m_1}\left(D_1\right)^{j_2-m_2}\left(D\right)^{2j_3-m_3}
 \cdot\cdot\cdot
 (\vphi)(t)\nonumber\\
 &&+
 \sum_{k_1+\ldots+k_N=N}
 (2\nu+1)^{k_1+k_3+\ldots}
 \sum_{s_1=0}^{k_1}\sum_{s_2=0}^{2k_2}\sum_{s_3=0}^{k_3}
 \cdot\cdot\cdot
 \binom{k_1}{s_1}\binom{2k_2}{s_2}\binom{k_3}{s_3}
 \cdot\cdot\cdot\nonumber\\
 &&\left(D_1\right)^{s_1}\left(D\right)^{s_2}\left(D_1\right)^{s_3}
 \cdot\cdot\cdot
 (\psi)(t)
 \cdot
 \left(D_1\right)^{k_1-s_1}\left(D\right)^{2k_2-s_2}\left(D_1\right)^{k_3-s_3}
 \cdot\cdot\cdot
 (\vphi)(t).
 \ena
 Next, for $\nu\in[-1/2,0),\, F(t)=f_\vep(t)g_\vep(t)$ by \eqref{E6.20aa},
 where $f_\vep \in E_{1-\vep,1,e}\cap L_\iy(\R^1)$ and
 $g_\vep \in E_{\vep,1,e}\cap L_\iy(\R^1)$.
 In addition,
 \beq\label{E6.20b2}
\left(Be_\nu\right)^N(F)(t)
=\left(Be_\nu\right)^N(f_\vep)(t)
+\left(Be_\nu\right)^N(f_\vep\cdot(g_\vep-1))(t).
\eeq
To estimate
\ba
\left\|\left(Be_\nu\right)^N
(f_\vep\cdot(g_\vep-1))\right\|_{L_\iy([-1,1])},
\ea
we use identity \eqref{E6.20b1} for $\vphi=f_\vep$
 and $\psi=g_\vep-1$. Then the uniform norm on
 $[-1,1]$ of all terms in \eqref{E6.20b1} with
 $m_1=\ldots=m_N=0$ and $s_1=\ldots=s_N=0$ can be estimated
 by Bernstein-type inequalities \eqref{E5.3a}
 and by the relations
 \ba
 \left\|g_\vep-1\right\|_{L_\iy([-1,1])}
 =\left\|g_\vep-g_\vep(0)\right\|_{L_\iy([-1,1])}
 \le C\vep.
 \ea
 All those estimates do not exceed
  $C\vep\left\|f\right\|_{L_\iy(\R^1)}$.
  The norms of other terms in \eqref{E6.20b1} do not exceed
  $C\vep\left\|f\right\|_{L_\iy(\R^1)}$ by Lemma \ref{L5.2a}
  (a) as well.
  Combining these estimates, we obtain
  \beq\label{E6.20b3}
  \eta_{\nu,N}\le \left\|\left(Be_\nu\right)^N
(f_\vep\cdot(g_\vep-1))\right\|_{L_\iy([-1,1])}
\le C_9\vep\left\|f\right\|_{L_\iy(\R^1)}.
\eeq
Therefore, setting $t=0$ in \eqref{E6.20b2},
we obtain \eqref{E6.20b} for $N\ge 1$ and $\nu\in[-1/2,0)$
from \eqref{E6.20b2} and \eqref{E6.20b3}.

Next, we use even polynomials $P_{2n}\in\PP_{2n,1,e}$
 from Lemma \ref{L5.2}
 such that for every
  $l\in\Z_+^1,$
  \bna\label{E6.20c}
   \lim_{n\to\iy} P_{2n}^{(l)}= F^{(l)},\qquad
 \lim_{n\to\iy} \left(Be_{\nu}\right)^l(P_{2n})= \left(Be_{\nu}\right)^l(F),
 \ena
 uniformly on each compact subset of $\R^1$.
 In addition, relation \eqref{E5.2a} of Lemma \ref{L5.2}
 for $N=0,\,r=p,\,\al=2\nu+1,$ and $\be=\nu$ shows that
 \beq\label{E6.20d}
 \lim_{n\to\iy}\|F-P_{2n}\|
 _{L_{p,\vert t\vert^{2\nu+1}
  \left(1-(t/(2\tau n))^2\right)^\nu}([-2\tau n,2\tau n])}=0.
  \eeq
 By Lemma \ref{L5.5} (a), there exists $R_n\in\PP_{n,1}$
such that for $b(n)=2\tau n$,
\ba
 P_{2n}(t)=P_{2n,b(n)}(t)=R_n\left(1-2(2\tau n)^{-2}t^2\right).
 \ea
 Then relations \eqref{E6.20c} show that we can use
 Lemma \ref{L5.5} (c) for $n_s=s,\,s\in\N$. Therefore, we obtain by \eqref{E5.15}
 \bna\label{E6.21}
&&\left\vert \left(Be_{\nu}\right)^N(F)(0)\right\vert
 =\lim_{n\to\iy}n^{-2N}\left\vert\left(Ge_{\nu+1/2}\right)^N(R_{n})(1)\right\vert\nonumber\\
 &&\le \liminf_{n\to\iy}n^{-2N}\left(
 \C_1\left(\PP_{n,1}
     ,\,L_{p,(1-u^2)^\nu}([-1,1]),\,
     \left(Ge_{\nu+1/2}\right)^N\right)
     \|R_n\|_{L_{p,(1-u^2)^\nu}([-1,1])}\right).
 \ena
Further, using the substitution $u=1-2(2\tau n)^{-2}t^2$,
"triangle" inequality \eqref{E2.7a}, and relation \eqref{E6.20d},
 we obtain for $p\in(0,\iy]$
and $\nu\in[-1/2,\iy)$
\bna\label{E6.22}
&&\limsup_{n\to\iy}\left(n^{(2\nu+2)/p}\|R_n\|_{L_{p,(1-u^2)^\nu}([-1,1])}\right)\nonumber\\
&&=2^{-1/p}\tau^{-(2\nu+2)/p}\limsup_{n\to\iy}
\left(\int_{-2\tau n}^{2\tau n}\left\vert P_{2n}(t)\right\vert^p\vert t\vert^{2\nu+1}\left(
1-\frac{t^2}{(2\tau n)^2}\right)^\nu dt\right)^{1/p}\nonumber\\
&&\le 2^{-1/p}\tau^{-(2\nu+2)/p}
\left(
\lim_{n\to\iy}
\left(\int_{-2\tau n}^{2\tau n}\left\vert P_{2n}(t)-F(t)\right\vert^p\vert t\vert^{2\nu+1}\left(
1-\frac{t^2}{(2\tau n)^2}\right)^\nu dt\right)^{\tilde{p}/p}\right.\nonumber\\
&&+\left.\limsup_{n\to\iy}
\left(\int_{-2\tau n}^{2\tau n}\left\vert F(t)\right\vert^p\vert t\vert^{2\nu+1}\left(
1-\frac{t^2}{(2\tau n)^2}\right)^\nu dt\right)^{\tilde{p}/p}\right)^{1/\tilde{p}}\nonumber\\
&&= 2^{-1/p}\tau^{-(2\nu+2)/p}\limsup_{n\to\iy}
\left(\int_{-2\tau n}^{2\tau n}\left\vert F(t)\right\vert^p\vert t\vert^{2\nu+1}\left(
1-\frac{t^2}{(2\tau n)^2}\right)^\nu dt\right)^{1/p}.
\ena
Next, we prove the estimate
\beq\label{E6.23}
\limsup_{n\to\iy}
\left(\int_{-2\tau n}^{2\tau n}\left\vert F(t)\right\vert^p\vert t\vert^{2\nu+1}\left(
1-\frac{t^2}{(2\tau n)^2}\right)^\nu dt\right)^{1/p}
\le (1-\vep)^{-(2\nu+2)/p}\|f\|_{L_{p,\vert t\vert^{2\nu+1}}(\R^1)}.
\eeq
It suffices to prove this inequality for $p\in(0,\iy)$.
For $\nu\ge 0$ inequality \eqref{E6.23} follows immediately from \eqref{E6.20a}.
If $\nu\in[-1/2,0)$,
then for any $\de\in(0,1)$ and $p\in(0,\iy)$,
\beq\label{E6.23a}
\int_{-2\tau n}^{2\tau n}\left\vert F(t)\right\vert^p\vert t\vert^{2\nu+1}\left(
1-\frac{t^2}{(2\tau n)^2}\right)^\nu dt
=2\left(\int_0^{2\de \tau n}+\int_{2\de \tau n}^{2\tau n}\right)
=2(I_1(n)+I_2(n)),
\eeq
where by \eqref{E6.20a},
\beq\label{E6.23b}
\limsup_{n\to\iy}I_1(n)
\le (1-\de^2)^{\nu}(1-\vep)^{-(2\nu+2)}
\int_0^{\iy}\vert f(t)\vert^p  t^{2\nu+1}\,dt,
\eeq
and by \eqref{E6.20aa},
\bna\label{E6.23c}
\limsup_{n\to\iy}I_2(n)
&=&\limsup_{n\to\iy}(2\tau n)^{2\nu+2}
\int_\de^1\left\vert F(2\tau n u)\right\vert^p u^{2\nu+1}
(1-u^2)^\nu du\nonumber\\
&\le& (1/2)B(\nu+1,\nu+1)\|f\|_{L_\iy(\R^1)}^p\limsup_{n\to\iy}(2\tau n)^{2\nu+2}
(2\vep \tau \de n/d)^{-d\,p}=0,
\ena
since $d\,p>2\nu+2$.

Collecting relations  \eqref{E6.23a}, \eqref{E6.23b}, and \eqref{E6.23c}, we obtain
\bna\label{E6.23d}
&&\limsup_{n\to\iy}
\left(\int_{-2\tau n}^{2\tau n}\left\vert F(t)\right\vert^p\vert t\vert^{2\nu+1}\left(
1-\frac{t^2}{(2\tau n)^2}\right)^\nu dt\right)^{1/p}\nonumber\\
&&\le (1-\de^2)^{\nu/p}(1-\vep)^{-(2\nu+2)/p}\|f\|_{L_{p,\vert t\vert^{2\nu+1}}(\R^1)}.
\ena
Letting $\de\to 0+$ in \eqref{E6.23d} we arrive at \eqref{E6.23}.

Combining \eqref{E6.20b} and \eqref{E6.21} with \eqref{E6.22} and \eqref{E6.23},
 we obtain for $p\in(0,\iy]$
and $\nu\in[-1/2,\iy)$
 \bna\label{E6.24}
 &&\left\vert(\Lambda_\nu(\vep))^{2N}\left( Be_{\nu}\right)^N(f)(0)
 +\eta_{\nu,N}(\vep)\right\vert
 \nonumber\\
     &&\le 2^{-1/p}\tau^{-(2\nu+2)/p}(1-\vep)^{-(2\nu+2)/p}\nonumber\\
      &&\times\liminf_{n\to\iy}n^{-2N-(2\nu+2)/p}
 \C_1\left(\PP_{n,1}
     ,\,L_{p,(1-u^2)^\nu}([-1,1]),\,
     \left(Ge_{\nu+1/2}\right)^N\right)
     \|f\|_{L_{p,\vert t\vert^{2\nu+1}}(\R^1)}\nonumber\\
     &&+C_9\vep\|f\|_{L_{\iy}(\R^1)}.
     \ena
 Letting $\tau\to 1-$ and $\vep\to 0+$ in \eqref{E6.24},
 we obtain
 \ba
 &&\left\vert\left( Be_{\nu}\right)^N(f)(0)\right\vert
 /\|f\|_{L_{p,\vert t\vert^{2\nu+1}}(\R^1)}\\
 &&\le  2^{-1/p}\liminf_{n\to\iy}n^{-2N-(2\nu+2)/p}
 \C_1\left(\PP_{n,1}
     ,\,L_{p,(1-u^2)^\nu}([-1,1]),\,
     \left(Ge_{\nu+1/2}\right)^N\right).
 \ea
 Hence we arrive at \eqref{E6.20} for $p\in(0,\iy]$
 and $\nu\in[-1/2,\iy)$.

 Further, we prove the inequality
\bna\label{E6.25}
&& \limsup_{n\to\iy} n^{-2N-(2\nu+2)/p}
     \C_1\left(\PP_{n,1}
     ,\,L_{p,(1-u^2)^\nu}([-1,1]),\,
     \left(Ge_{\nu+1/2}\right)^N\right)\nonumber\\
     &&\le
     2^{1/p}\C_0\left(E_{1,1,e}\cap L_{p,\vert t\vert^{2\nu+1}}(\R^1),
     \,L_{p,\vert t\vert^{2\nu+1}}(\R^1),\,
     \left(Be_{\nu}\right)^N\right)
     \ena
by constructing a nontrivial function $f_0\in E_{1,1,e}
  \cap L_{p,\vert t\vert^{2\nu+1}}(\R^1)$
   such that
 \bna \label{E6.26}
&&\limsup_{n\to\iy} n^{-2N-(2\nu+2)/p}
     \C_1\left(\PP_{n,1}
     ,\,L_{p,(1-u^2)^\nu}([-1,1]),\,
     \left(Ge_{\nu+1/2}\right)^N\right)\nonumber\\
  &&\le 2^{1/p}\left\vert \left(Be_{\nu}\right)^N(f_0)(0)\right\vert/
\|f_0\|_{L_{p,\vert t\vert^{2\nu+1}}(\R^1)}\nonumber\\
  &&\le 2^{1/p}\C_0\left(E_{1,1,e}\cap L_{p,\vert t\vert^{2\nu+1}}(\R^1),
     \,L_{p,\vert t\vert^{2\nu+1}}(\R^1),\,
     \left(Be_{\nu}\right)^N\right).
 \ena
 Then inequalities \eqref{E6.20} and \eqref{E6.25} imply  \eqref{E4.5}.
 In addition, $f_0$ is an extremal function in \eqref{E4.5},
that is, relation \eqref{E4.5a} is valid.

It remains to construct a nontrivial function $f_0$,
satisfying \eqref{E6.26}.
We first note that
\beq \label{E6.27}
\inf_{n\ge 2N+2}n^{-2N-(2\nu+2)/p}
\C_1\left(\PP_{n,1}
     ,\,L_{p,(1-u^2)^\nu}([-1,1]),\,
     \left(Ge_{\nu+1/2}\right)^N\right)
     \ge C_{10}(p,N,\nu).
\eeq
This inequality follows immediately from  \eqref{E6.20}.
Let $R_n\in\PP_{n, 1}$
 be a polynomial, satisfying the equality
\beq \label{E6.28}
\C_1\left(\PP_{n,1}
     ,\,L_{p,(1-u^2)^\nu}([-1,1]),\,
     \left(Ge_{\nu+1/2}\right)^N\right)
=\frac{\left\vert \left(Ge_{\nu+1/2}\right)^N(R_n)(1)\right\vert}
{\|R_n\|_{L_{p,(1-u^2)^\nu}([-1,1])}},\quad n\in\N.
\eeq
The existence of an extremal polynomial $R_n$ in \eqref{E6.28}
can be proved by the standard compactness argument.
We can assume that
\beq \label{E6.29}
 n^{-2N} \left(Ge_{\nu+1/2}\right)^N(R_n)(1)=1,\qquad n\in\N.
\eeq
Next, setting $Q_{2n}(t):=R_n(1-2(2n)^{-2}t^2)$, we have from
\eqref{E6.28}, \eqref{E6.29}, and \eqref{E6.27} that for $n\ge 2N+2$,
\bna \label{E6.30}
&&\|Q_{2n}\|_{L_{p,\vert t\vert^{2\nu+1}(1-(t/(2n))^2)^\nu}([-2n,2n])}
=2^{1/p}n^{(2\nu+2)/p}
\|R_{n}\|_{L_{p,(1-u^2)^\nu}([-1,1])}\nonumber\\
&&=2^{1/p}\left(n^{-2N-(2\nu+2)/p}
\C_1\left(\PP_{n,1}
     ,\,L_{p,(1-u^2)^\nu}([-1,1]),\,
     \left(Ge_{\nu+1/2}\right)^N\right)\right)^{-1}\nonumber\\
     &&\le  2^{1/p}/C_{10}(p,N,\nu).
     \ena
Further, $Q_{2n}\in\PP_{2n,1}$, and combining inequality \eqref{E5.6}
of Lemma \ref{L5.4}
for $\al=2\nu+1$ and $\be=\nu$ with
 \eqref{E6.29} we obtain
for any $\vep\in(0,1/2)$ and any $k\in\Z_+^1$,
\bna \label{E6.31}
\left\vert Q_{2n}^{(k)}(0)\right\vert
&\le& C_4(2\nu+1,\nu,p,\vep)(1-\vep)^{-k }
\|Q_{2n}\|_{L_{p,\vert t\vert^{2\nu+1}(1-(t/(2n))^2)^\nu}([-2n,2n])}\nonumber\\
&\le& \left(2^{1/p}C_4(2\nu+1,\nu,p,\vep)/C_{10}(p,N,\nu)\right)(1-\vep)^{-k}.
\ena
Let $\{n_r\}_{r=1}^\iy$ be a subsequence of natural numbers such that
\bna \label{E6.32}
&&\limsup_{n\to\iy}\,
n^{-2N-(2\nu+2)/p}
\C_1\left(\PP_{n,1}
     ,\,L_{p,(1-u^2)^\nu}([-1,1]),\,
     \left(Ge_{\nu+1/2}\right)^N\right)\nonumber\\
     &&=\lim_{r\to\iy}
     (n_r)^{-2N-(2\nu+2)/p}
     \C_1\left(\PP_{n_r,1}
     ,\,L_{p,(1-u^2)^\nu}([-1,1]),\,
     \left(Ge_{\nu+1/2}\right)^N\right).
\ena
Inequality \eqref{E6.31} shows that the polynomial sequence
$\{Q_{2n_r}\}_{r=1}^\iy\subseteq \E_1$ satisfies the conditions
of Lemma \ref{L5.3}. Therefore, there exist
 a function $f_0\in E_{1,1,e}$ and a subsequence
 $\{Q_{2n_{r_s}}\}_{s=1}^\iy$ such that for all $l\in\Z_+^1$,
 \beq \label{E6.33}
 \lim_{s\to\iy}Q_{2n_{r_s}}^{(l)}(t)=f_0^{(l)}(t),
 \qquad \lim_{s\to\iy}(Be_\nu)^l\left(Q_{2n_{r_s}}\right)(t)
 =(Be_\nu)^l(f_0)(t),
 \eeq
 uniformly on any interval $[-A,A],\,A>0$.
 In addition, it follows from \eqref{E6.33} for $l=0$ that
 for $p\in(0,\iy],$ and $\nu\in[-1/2,\iy)$,
 \bna\label{E6.33a}
 &&\lim_{s\to\iy}
\|f_0-Q_{2n_{r_s}}\|_{L_{p,\vert t\vert^{2\nu+1}
\left(1-(t/(2n_{r_s})\right)^2)^\nu}([-A,A])}\nonumber\\
&&\le \lim_{s\to\iy}\|f_0-Q_{2n_{r_s}}\|_{L_\iy([-A,A])}
\limsup_{s\to\iy}\left(
\int_{-A}^{A}\vert t\vert^{2\nu+1}\left(
1-\frac{t^2}{(2 n_{r_s})^2}\right)^\nu dt\right)^{1/p}
=0.
\ena
 Then relations \eqref{E6.33}
 show that we can use
 Lemma \ref{L5.5} (c) for $n_s=n_{r_s},\,s\in\N$.
 Therefore, we obtain by \eqref{E5.15}
 \beq \label{E6.34}
  \lim_{s\to\iy}n_{r_s}^{-2N} \left(Ge_{\nu+1/2}\right)^N(R_{n_{r_s}})(1)
  =\left(Be_{\nu}\right)^N(f_0)(0).
  \eeq
  It follows from \eqref{E6.29} and \eqref{E6.34} that
  \beq \label{E6.35}
  \left(Be_{\nu}\right)^N(f_0)(0)=1.
 \eeq
In addition,   applying "triangle" inequality \eqref{E2.7a}
and using \eqref{E6.33a} and \eqref{E6.30},
we obtain for any interval $[-A,A],\,A>0$,
\bna \label{E6.36}
&&\|f_0\|_{L_{p,\vert t\vert^{2\nu+1}}([-A,A])}^{\tilde{p}}
=\lim_{s\to\iy}\|f_0\|_{L_{p,\vert t
\vert^{2\nu+1}\left(1-(t/(2n_{r_s}))^2\right)^\nu}([-A,A])}^{\tilde{p}}\nonumber\\
&&\le \limsup_{s\to\iy}
\left(\|f_0-Q_{2n_{r_s}}\|_{L_{p,\vert t\vert^{2\nu+1}
\left(1-(t/(2n_{r_s}))^2\right)^\nu}([-A,A])}^{\tilde{p}}
+\|Q_{2n_{r_s}}\|_{L_{p,\vert t\vert^{2\nu+1}
\left(1-(t/(2n_{r_s}))^2\right)^\nu}([-A,A])}^{\tilde{p}}\right)\nonumber\\
&&\le \limsup_{s\to\iy}\|Q_{2n_{r_s}}
\|_{L_{p,\vert t\vert^{2\nu+1}
\left(1-(t/(2n_{r_s}))^2\right)^\nu}([-2n_{r_s},2n_{r_s}])}^{\tilde{p}}
\le \left(2^{1/p}/C_{10}(p,N,\nu)\right)^{\tilde{p}}.
\ena
Therefore, $f_0$ is a nontrivial function from
$E_{1,1.e}\cap L_{p,\vert t\vert^{2\nu+1}}(\R^1)$,
by \eqref{E6.36} and \eqref{E6.35}. Thus for any interval
$[-A,A],\,A>0$, we obtain from \eqref{E6.30}, \eqref{E6.32},
 \eqref{E6.33a},
 and \eqref{E6.35}
\bna \label{E6.37}
&&2^{-1/p}\limsup_{n\to\iy} n^{-2N-(2\nu+2)/p}
\C_1\left(\PP_{n,1}
     ,\,L_{p,(1-u^2)^\nu}([-1,1]),\,
     \left(Ge_{\nu+1/2}\right)^N\right)
     \nonumber\\
&&=\lim_{s\to\iy}\left(\|Q_{2n_{r_s}}\|
_{L_{p,\vert t\vert^{2\nu+1}\left(1-(t/(2n_{r_s}))^2\right)^\nu}
([-2n_{r_s},2n_{r_s}])}\right)^{-1}\nonumber\\
&&\le \lim_{s\to\iy}\left(\|Q_{n_{r_s}}\|
_{L_{p,\vert t\vert^{2\nu+1}
 \left(1-(t/(2n_{r_s}))^2\right)^\nu}([-A,A])}\right)^{-1}\nonumber\\
&&=\left\vert \left(Be_{\nu}\right)^N (f_0)(0)\right\vert/
\|f_0\|_{L_{p,\vert t\vert^{2\nu+1}}([-A,A])}.
\ena
Finally, letting $A\to \iy$ in \eqref{E6.37}, we arrive at \eqref{E6.26}
for $\nu\ge -1/2,\,N\in\Z_+^1$, and $p\in(0,\iy]$.
\hfill $\Box$
\vspace{.12in}\\
\textbf{Acknowledgement.} We are grateful to both anonymous referees
 for valuable suggestions.

\end{document}